\documentclass{amsart}
\usepackage[all]{xy}
\usepackage{tikz}
 \tikzstyle{int}=[circle, draw,fill=black,outer sep=0,minimum size=3pt, inner sep=0]
  \tikzstyle{ext}=[circle, draw=black,outer sep=0,inner sep=1pt]
\usepackage{latexsym}
\usepackage{amssymb}
\usepackage{amsmath}
\usepackage{amsthm}
\usepackage{amscd}
\usepackage{fancyhdr}
\usepackage{fullpage}
\usepackage{mathrsfs}
\input cyracc.def
\xyoption{arc}

\def\id{{\mbox{1 \hskip -8pt 1}}}
\newcommand{\sgn}{{\mathit s  \mathit g\mathit  n}}
 \newcommand{\lon}{\longrightarrow}
 \newcommand{\bu}{\bullet}
 
 \newcommand{\rar}{\rightarrow}

 \newcommand{\Z}{{\mathbb Z}}
 \newcommand{\bS}{{\mathbb S}}

 \newcommand{\K}{{\mathbb K}}

\newcommand{\GC}{\mathsf{GC}}

 \newcommand{\ot}{\otimes}

\newcommand{\sC}{{\mathsf C}}

\newcommand{\sG}{{\mathsf G}}
\newcommand{\sP}{{\mathsf P}}

\newcommand{\Def}{{\mathsf D\mathsf e\mathsf f }}

 \newcommand{\Beq}{\begin{equation}}
 \newcommand{\Eeq}{\end{equation}}
 \newcommand{\Beqr}{\begin{eqnarray}}
 \newcommand{\Eeqr}{\end{eqnarray}}
 \newcommand{\Beqrn}{\begin{eqnarray*}}
 \newcommand{\Eeqrn}{\end{eqnarray*}}
 \newcommand{\Ba}{\begin{array}}
 \newcommand{\Ea}{\end{array}}
 \newcommand{\Bi}{\begin{itemize}}
 \newcommand{\Ei}{\end{itemize}}
 \newcommand{\Bc}{\begin{center}}
 \newcommand{\Ec}{\end{center}}
 
 \newcommand{\fg}{{\mathfrak g}}

\newcommand{\ft}{{\mathfrak t}}
\newcommand{\fr}{{\mathfrak r}}

 \newcommand{\cC}{{\mathcal C}}
 \newcommand{\caD}{{\mathcal D}}
 \newcommand{\cE}{{\mathcal E}}
 \newcommand{\cF}{{\mathcal F}}
 \newcommand{\cG}{{\mathcal G}}
 \newcommand{\caH}{{\mathcal H}}

 \newcommand{\caL}{{\mathcal L}}
 \newcommand{\cM}{{\mathcal M}}
 
 \newcommand{\cP}{{\mathcal P}}
 
 \newcommand{\cR}{{\mathcal R}}
 
 \newcommand{\cT}{{\mathcal T}}

 \newcommand{\Ga}{\Gamma}

 \newcommand{\la}{\lambda}

 \newcommand{\Hom}{{\mathrm H\mathrm o\mathrm m}}
 
 \newcommand{\sip}{\smallskip}
 \newcommand{\bip}{\bigskip}
 \newcommand{\mip}{\vspace{2.5mm}}

\newcommand{\hoe}{\mathrm{hoe}}
\newcommand{\Aut}{\mathrm{Aut}}
 \newcommand{\GCor}{\GC^{or}}
 \newcommand{\hGCor}{\widehat{\GC}^{or}}

 \newcommand{\LieBi}{{\caL \mathit{ieb}}}
  \newcommand{\hLieBi}{\widehat{\LieBi}}
  \newcommand{\hLoB}{\widehat{\LieBi}^{\Ba{c}\vspace{-1mm}_{\hspace{-2mm}\diamond} \Ea}}

 \newcommand{\LieBiP}{\LieBi P}

  \newcommand{\LB}{\mathcal{L}\mathit{ieb}}
\newcommand{\LoB}{\mathcal{L}\mathit{ieb}^\diamond}

\newcommand{\LBcd}{\mathcal{L}\mathit{ieb}_{c,d}}

\newcommand{\wLBcd}{\widehat{\mathcal{L}\mathit{ieb}}_{c,d}}

\newcommand{\HoLBcd}{\mathcal{H}\mathit{olieb}_{c,d}}

\newcommand{\wHoLB}{\widehat{\mathcal{H}\mathit{olieb}}}
\newcommand{\wHoLBcd}{\widehat{\mathcal{H}\mathit{olieb}}_{c,d}}

\newcommand{\LoBcd}{\mathcal{L}\mathit{ieb}_{c,d}^\diamond}

\newcommand{\wLoBcd}{\widehat{\mathcal{L}\mathit{ieb}}_{c,d}^{\Ba{c}\vspace{-1mm}_{\hspace{-2mm}\diamond} \Ea}}
\newcommand{\HoLoBd}{\mathcal{H}\mathit{olieb}_{d}^\diamond}
\newcommand{\HoLoBcd}{\mathcal{H}\mathit{olieb}_{c,d}^\diamond}
\newcommand{\HoLoB}{\mathcal{H}\mathit{olieb}^\diamond}
\newcommand{\HoLB}{\mathcal{H}\mathit{olieb}}

\newcommand{\wHoLoB}{\widehat{\mathcal{H}\mathit{olieb}}^{\Ba{c}\vspace{-1mm}_{\hspace{-2mm}\diamond} \Ea}}
\newcommand{\wHoLoBcd}{\widehat{\mathcal{H}\mathit{olieb}}_{c,d}^{\Ba{c}\vspace{-1mm}_{\hspace{-2mm}\diamond} \Ea}}

 \newcommand{\grt}{\mathfrak{grt}}
 \newcommand{\Der}{\mathrm{Der}}

 \newcommand{\gr}{\mathrm{gr}}
 
\theoremstyle{plain}
\swapnumbers

\newtheorem{prop-def}[theorem]{Proposition-definition}

\newtheorem{main-theorem}{Main~Theorem}[section]
\newtheorem{section-theorem}{Theorem}[section]
\newtheorem{section-corollary}{Corollary}[section]

\theoremstyle{definition}

\begin{document}

\sloppy

 \newenvironment{proo}{\begin{trivlist} \item{\sc {Proof.}}}
  {\hfill $\square$ \end{trivlist}}

\long\def\symbolfootnote[#1]#2{\begingroup%
\def\thefootnote{\fnsymbol{footnote}}\footnote[#1]{#2}\endgroup}

  \title{Deformation theory of  Lie bialgebra properads}

\author{Sergei~Merkulov}
\address{Sergei~Merkulov:  Mathematics Research Unit, Luxembourg University,  Grand Duchy of Luxembourg }
\email{sergei.merkulov@uni.lu}

\author{Thomas~Willwacher}
\address{Thomas~Willwacher: Institute of Mathematics, University of Zurich, Zurich, Switzerland}
\email{thomas.willwacher@math.uzh.ch}



 \begin{abstract}
 We compute the homotopy derivations of the properads governing even and odd Lie bialgebras as well as involutive Lie bialgebras.
 The answer may be expressed in terms of the Kontsevich graph complexes.
 In particular, this shows that the Grothendieck-Teichm\"uller group acts faithfully (and essentially transitively) on the  completions of the properads governing even Lie bialgebras and involutive Lie bialgebras, up to homotopy. This shows also that by contrast to the even case the properad governing odd Lie bialgebras admits precisely one non-trivial automorphism --- the standard rescaling automorphism,
 and that it has precisely one non-trivial deformation which we describe explicitly.
\end{abstract}
 \maketitle

{\large
\section{\bf Introduction}
}
\label{sec:introduction}

\subsection{Deformation theory of Lie bialgebras and graph complexes} An {\em even}\, Lie bialgebra is a vector space which carries both a Lie algebra and a Lie coalgebra structure of the same $\Z_2$ parity, satisfying a certain compatibility relation. Lie bialgebras were introduced by Drinfeld in \cite{D1}  in the context of the
theory of Yang-Baxter equations. They have since seen numerous application, in particular in the theory
of Hopf algebra deformations of universal enveloping algebras, cf. \cite{ES} and
references therein.

 \sip

If the composition of the cobracket and bracket of a Lie bialgebra is zero, the Lie bialgebra is called \emph{involutive}.
This additional condition is satisfied in many interesting examples studied in homological algebra, string topology, symplectic field theory,  Lagrangian Floer theory of higher genus, and the theory of cohomology groups $H(\cM_{g,n})$ of moduli spaces of algebraic curves with labelings of punctures  {\em skewsymmetrized} \cite{Ch,CFL,MW,Tu, Sch}.

\sip

{\em Odd}\, Lie bialgebras are, by definition, the ones in which Lie brackets and cobrackets have opposite $\Z_2$
parities.  They have seen applications in Poisson geometry \cite{Me1}, deformation quantization
of Poisson structures
\cite{Me2} and, most surprisingly, in the theory of cohomology groups $H(\cM_{g,n})$ of moduli spaces of algebraic curves with labelings of punctures {\em symmetrized}\, \cite{MW}.

\sip

We study in this paper the deformation theory of Lie bialgebras in both the even and odd cases.
Let $\LB_{c,d}$ denote the properad governing Lie bialgebras with Lie brackets in degree $1-c$
and Lie cobrackets in degree $1-d$ so that the case $c+d\in 2\Z$ corresponds to even Lie bialgebras, and the case $c+d\in 2\Z+1$ to odd ones. The involutivity condition is non-trivial only in the even case, and we denote by $\LoB_{c,d}$ (with $c+d\in 2\Z$ by default)   the properad governing involutive Lie bialgebras.

\sip

The main purpose of this note is to finish the study of the homotopy derivations of the properads $\LB_{c,d}$ and $\LoB_{c,d}$ initiated in \cite{CMW} where a minimal resolution $\HoLoB_{c,d}$ of
$\LoB_{c,d}$ was constructed; minimal resolutions $\HoLB_{c,d}$ of $\LB_{c,d}$ have been found
earlier in  \cite{Ko,MaVo} for $c+d\in 2\Z$ and in \cite{Me1,Me2} for $c+d\in 2\Z +1$.
Informally speaking, one question we want to answer is which universal deformations of any  (involutive) Lie bialgebra structure one can construct, using only the bracket and cobracket operations.

\sip

The answers may be formulated in terms of the cohomology of the oriented graph complex $\GCor_d$,
whose elements are linear combinations of isomorphism classes of directed  graphs with no closed paths of directed edges. 
More concretely, in \cite{CMW} the authors describe maps of dg Lie algebras
$$
F \colon \GC_{c+d+1}^{or}\to \Der(\wHoLBcd) \ \ \ \ , \ \ \ \
F^\diamond \colon \GC_{c+d+1}^{or}[[\hbar]]\to \Der(\wHoLoBcd),
$$
where the symbol\  $\widehat{\ \ \ \ \ \ }$\  stands for the genus completion, $\hbar$ for a formal variable of homological degree $c+d$, and $\Der(\cP)$ for the dg Lie algebra of (continuous) derivations of a dg properad $\cP$.

\sip

 Precise definitions of all the properads and the graph complexes and in particular the dg Lie algebra structures on  $\GCor_{c+d+1}$ will be recalled in \S 2 below. Explicit formulae for the morphisms $F$ and $F^\diamond$ as well for the dg Lie algebra structures in $\Der(\wHoLBcd)$ and $\Der(\wHoLoBcd)$ are given in \S 3.

\subsection{Main Theorems}
The main results of the present note are the following two Theorems which are proven in \S 4.

\subsubsection{\bf Theorem} \label{thm:Fqiso}
{\em For any $c,d\in \Z$ the map
$$
F \colon \GC_{c+d+1}^{or}\to \Der(\wHoLBcd)
\vspace{-4mm}
$$
is a quasi-isomorphism, up to one class in  $\Der(\wHoLBcd)$
 represented by the series
 $\displaystyle
  \sum_{m,n}(m+n-2)
  \overbrace{
  \underbrace{
 \Ba{c}\resizebox{6mm}{!}  {\xy
(0,4.5)*+{...},
(0,-4.5)*+{...},
(0,0)*{\bu}="o",
(-5,5)*{}="1",
(-3,5)*{}="2",
(3,5)*{}="3",
(5,5)*{}="4",
(-3,-5)*{}="5",
(3,-5)*{}="6",
(5,-5)*{}="7",
(-5,-5)*{}="8",
\ar @{-} "o";"1" <0pt>
\ar @{-} "o";"2" <0pt>
\ar @{-} "o";"3" <0pt>
\ar @{-} "o";"4" <0pt>
\ar @{-} "o";"5" <0pt>
\ar @{-} "o";"6" <0pt>
\ar @{-} "o";"7" <0pt>
\ar @{-} "o";"8" <0pt>
\endxy}\Ea
 }_{n\times}
 }^{m\times}.
 $
}

\subsubsection{\bf Theorem}\label{thm:Fhbarqiso}
{\em  For any $c,d\in \Z$  with $c+d\in 2\Z$ the map
 $$
 F^\diamond \colon \GC_{c+d+1}^{or}[[\hbar]]\to \Der(\wHoLoBcd)
 $$
   is a quasi-isomorphisms, up to classes $T\K[[\hbar]]\in \Der(\wHoLoBcd)$, generated over $\K[[\hbar]]$ by the series
 \[
  T=
  \sum_{m,n,p}(m+n+2p-2) \hbar^{p}
   \overbrace{
 \underbrace{
 \Ba{c}\resizebox{8mm}{!}  { \xy
(0,4.5)*+{...},
(0,-4.5)*+{...},
(0,0)*+{_p}*\cir{}="o",
(-5,5)*{}="1",
(-3,5)*{}="2",
(3,5)*{}="3",
(5,5)*{}="4",
(-3,-5)*{}="5",
(3,-5)*{}="6",
(5,-5)*{}="7",
(-5,-5)*{}="8",
\ar @{-} "o";"1" <0pt>
\ar @{-} "o";"2" <0pt>
\ar @{-} "o";"3" <0pt>
\ar @{-} "o";"4" <0pt>
\ar @{-} "o";"5" <0pt>
\ar @{-} "o";"6" <0pt>
\ar @{-} "o";"7" <0pt>
\ar @{-} "o";"8" <0pt>
\endxy}\Ea
 }_{n\times}
 }^{m\times}.
 \]
}

\subsection{Some applications} The most important cases for applications correspond to (i) the case $c=d=1$ where one deals
with the classical (i.e.\ the ones in which with all the generators have homological degree zero)  Lie bialgebra properads $\LB:=\LB_{1,1}$ and $\LoB:=\LoB_{1,1}$, and the case (ii) $c=0$, $d=1$  where one gets a properad $\LB_{odd}:= \LB_{0,1}$ of odd Lie bialgebras which has the property \cite{Me1,Me2} that the representations of its minimal resolution $\HoLB_{odd}$ in a graded
vector space $V$ are in 1-1 correspondence with formal graded Poisson structures $\pi \in \cT_{poly}^{\geq 1}(V)$ on $V$
which vanish at $0\in V$.

\sip

As has been shown in \cite{Wi1,Wi2} the cohomology of  the oriented graph complexes $H(\GCor_{d+1})$ may be identified with the cohomology of the ``plain" graph complexes\footnote{The superscript $2$ in the symbol $\GC_d^2$ means that we consider graphs with at least bivalent vertices, while the symbol $\GC_d$ is reserved traditionally  for a complex of graphs with at least trivalent vertices.} $\GC_{d}^2$ introduced
(for $d=2$) by Kontsevich  in \cite{Kon}. The latter complexes have been studied in \cite{Wi1} where it was proven, in particular, that there is an isomorphism
of Lie algebras,
$$
H^0(\GC_{2}^2)= \fg\fr\ft_1,
$$
where $\fg\fr\ft_1$ is the Lie algebra of the Grothendieck-Teichm\"uller group $GRT_1$ introduced by Drinfeld in
\cite{D2}. Then, according to  \cite{Wi1,Wi2} and, respectively, \cite{CMW} , one concludes that
$$
H^0(\GCor_{3})=\fg\fr\ft_1  \ \ \mbox{and}\ \ \  H^0(\sG\sC_3^{or}[[\hbar]])=\fg\fr\ft_1
$$
and hence obtains (as corollaries to the Main Theorems {\ref{thm:Fqiso}} and {\ref{thm:Fhbarqiso}} above) the faithful  actions of the group $GRT_1$ on the completed properads $\wHoLB$ and $\wHoLoB$ and
and hence on their representations (the precise meaning of representation of a {\em completed}\, properad is given in \S {\ref{3: subsect on GCor_3 to Der(LieB)}} below).

\sip

It is a folklore conjecture that the cohomology group $H^1(\GC_{2}^2)=H^1(\GCor_{3})$ vanishes; if proven, the  Main Theorems would imply that the properads $\HoLB$ and $\HoLoB$ and of their completions
are {\em rigid}, i.e.\ admit no nontrivial deformations.

\sip

By contrast, one has $H^0(\GCor_2)=0$, $H^1(\GCor_2)=\K$ and $H^2(\GCor_2)=\K$. Then the  Main Theorem {\ref{thm:Fqiso}} says
that the properad  $\HoLB_{odd}$ and its genus completion admits precisely one homotopy non-trivial automorphism --  the standard rescaling automorphism, and that its completed
version $\widehat{\HoLB}_{odd}$ has  precisely  {one non-trivial} deformation which we describe explicitly in \S {{\ref{4: subsec on applications}}. This unique deformed version of $\widehat{\HoLB}_{odd}$ leads to the notion of {\em quantizable Poisson structures} \cite{MW3, KMW}, the ones whose deformation quantization is a trivial procedure not requiring a choice of an associator as in the case of deformation quantization of ordinary Poisson structures.

\subsection{Some notation}
 The set $\{1,2, \ldots, n\}$ is abbreviated to $[n]$;  its group of automorphisms is
denoted by $\bS_n$;
the trivial one-dimensional representation of
 $\bS_n$ is denoted by $\id_n$, while its one dimensional sign representation is
 denoted by $\sgn_n$.
The cardinality of a finite set $A$ is denoted by $\# A$.

\sip

We work throughout in the category of $\Z$-graded vector spaces over a field $\K$
of characteristic zero.
If $V=\oplus_{i\in \Z} V^i$ is a graded vector space, then
$V[k]$ stands for the graded vector space with $V[k]^i:=V^{i+k}$ and
and $s^k$ for the associated isomorphism $V\rar V[k]$; for $v\in V^i$ we set $|v|:=i$.
For a pair of graded vector spaces $V_1$ and $V_2$, the symbol $\Hom_i(V_1,V_2)$ stands
for the space of homogeneous linear maps of degree $i$, and
$\Hom(V_1,V_2):=\bigoplus_{i\in \Z}\Hom_i(V_1,V_2)$; for example, $s^k\in \Hom_{-k}(V,V[k])$.

\sip

For a
properad $\cP$ we denote by $\cP\{k\}$ the properad which is uniquely defined by
 the following property:
for any graded vector space $V$ a representation
of $\cP\{k\}$ in $V$ is identical to a representation of  $\cP$ in $V[k]$.
 The degree shifted operad of Lie algebras $\caL \mathit{ie}\{d\}$  is denoted by $\caL ie_{d+1}$.
  while its minimal resolution by $\caH \mathit{olie}_{d+1}$; representations of $\caL ie_{d+1}$ are vector spaces equipped with Lie brackets of degree $-d$.

\sip

For a right (resp., left) module $V$ over a group $G$ we denote by $V_G$ (resp.\
$_G\hspace{-0.5mm}V$)
 the $\K$-vector space of coinvariants:
$V/\{g(v) - v\ |\ v\in V, g\in G\}$ and by $V^G$ (resp.\ $^GV$) the subspace
of invariants: $\{\forall g\in  G\ :\  g(v)=v,\ v\in V\}$. If $G$ is finite, then these
spaces are canonically isomorphic as $char(\K)=0$.

\sip

For a vector space $V$ and a formal parameter $\hbar$ we denote by $V[[\hbar]]$ the topological vector space
of formal power series in $\hbar$ with coefficients in $V$, and by $\hbar^k V[[\hbar]]$ a vector space
of $V[[\hbar]]$ spanned by series of the form $\hbar^k f$ for some $f\in V[[\hbar]]$.
\mip

\subsection{Remark} A part of this paper contains an extended version of an appendix in the preprint
\cite{CMW} which  was removed from the final version of  loc.\ cit.\ following the recommendation of a referee.

\subsection*{Acknowledgements}
S.M. acknowledges support
by the Swedish Vetenskapr\aa det, grant 2012-5478.
T.W. acknowledges partial support by the Swiss National Science Foundation, grant 200021\_150012, and by the NCCR SwissMAP of the Swiss National Science Foundation.

\bip

{
\Large
\section{\bf Properads of  Lie bialgebras and graph complexes}\label{sec:preliminaries}
}

\sip

\subsection{Lie $n$-bialgebras} A  {\em Lie n-bialgebra}\, is a graded vector space $V$
equipped with linear maps,
$$
\vartriangle: V\rightarrow V\wedge V \ \ \ \mbox{and}\ \ \  [\ , \ ]: \wedge^2 (V[n])
\rightarrow V[n],
$$
such that the first operation $\vartriangle$ makes $V$ into a Lie coalgebra, the second operation
$[\ ,\ ]$ makes  $V[n]$ into a Lie algebra, and  the compatibility condition
$$
\vartriangle [a, b] = \sum a_1\otimes [a_2, b] +  [a,
b_1]\otimes b_2 - (-1)^{(|a|+n)(|b|+n)}( [b, a_1]\otimes a_2
+ b_1\otimes [b_2, a]),
$$
holds for any $a,b\in V$ with $\vartriangle a=:\sum a_1\otimes a_2$, $\vartriangle b=:\sum
b_1\otimes b_2$.
The case  $n=0$  gives us the ordinary definition of Lie bialgebra \cite{D1}.
The case $n=1$  is if interest because minimal resolutions of Lie 1-bialgebras control
local Poisson geometry \cite{Me1,Me2}. For $n$ even
it makes sense to introduce an {\em involutive Lie
$n$-bialgebra}\, as a Lie $n$-bialgebra $(V, [\ ,\ ], \vartriangle)$ such that the
composition map
$$
\Ba{ccccc}
V & \stackrel{\vartriangle}{\lon} & \Lambda^2V & \stackrel{[\ ,\ ]}{\lon} & V[-n]\\
a & \lon &    \sum a_1\otimes a_2 &\lon & [a_1,a_2]
\Ea
$$
vanishes (for odd $n$ this condition is trivial for symmetry reasons).

\sip

\subsection{Properads of (involutive) Lie bialgebras.} Let $\caL ieb_{n}$
(resp.\ $\caL ieb^\diamond_{n}$) denote the properad  of (resp.\ involutive) Lie $n$-bialgebras.
Let us consider their degree shifted versions,
\[
 \caL ieb_{c,d} = \caL ieb_{c+d-2}\{1-c\}, \ \ \ \  \caL ieb_{c,d}^\diamond = \caL ieb_{c+d-2}\{1-c\}
\]
in which the cobracket generator has degree $1-c$, and the bracket generator degree $1-d$. It is worth emphasizing that
the symbol $\caL ieb_{c,d}^\diamond$ tacitly assumes that $c+d\in 2\Z$, i.e.\ that the numbers $c$ and $d$ have the same parity. Let us describe these properads and their minimal resolutions explicitly.

\sip

By definition,  $\LBcd$ is a quadratic properad given as the quotient,
$$
\LB_{c,d}:=\cF ree\langle E\rangle/\langle\cR\rangle,
$$
of the free properad generated by an  $\bS$-bimodule $E=\{E(m,n)\}_{m,n\geq 1}$ with
 all $E(m,n)=0$ except
$$
E(2,1):=\id_1\ot \sgn_2^{c}[c-1]=\mbox{span}\left\langle
\Ba{c}\begin{xy}
 <0mm,-0.55mm>*{};<0mm,-2.5mm>*{}**@{-},
 <0.5mm,0.5mm>*{};<2.2mm,2.2mm>*{}**@{-},
 <-0.48mm,0.48mm>*{};<-2.2mm,2.2mm>*{}**@{-},
 <0mm,0mm>*{\circ};<0mm,0mm>*{}**@{},
 <0.5mm,0.5mm>*{};<2.7mm,2.8mm>*{^{_2}}**@{},
 <-0.48mm,0.48mm>*{};<-2.7mm,2.8mm>*{^{_1}}**@{},
 \end{xy}\Ea
=(-1)^{c}
\Ba{c}\begin{xy}
 <0mm,-0.55mm>*{};<0mm,-2.5mm>*{}**@{-},
 <0.5mm,0.5mm>*{};<2.2mm,2.2mm>*{}**@{-},
 <-0.48mm,0.48mm>*{};<-2.2mm,2.2mm>*{}**@{-},
 <0mm,0mm>*{\circ};<0mm,0mm>*{}**@{},
 <0.5mm,0.5mm>*{};<2.7mm,2.8mm>*{^{_1}}**@{},
 <-0.48mm,0.48mm>*{};<-2.7mm,2.8mm>*{^{_2}}**@{},
 \end{xy}\Ea
   \right\rangle
$$
$$
E(1,2):= \sgn_2^{d}\ot \id_1[d-1]=\mbox{span}\left\langle
\Ba{c}\begin{xy}
 <0mm,0.66mm>*{};<0mm,3mm>*{}**@{-},
 <0.39mm,-0.39mm>*{};<2.2mm,-2.2mm>*{}**@{-},
 <-0.35mm,-0.35mm>*{};<-2.2mm,-2.2mm>*{}**@{-},
 <0mm,0mm>*{\circ};<0mm,0mm>*{}**@{},
   <0.39mm,-0.39mm>*{};<2.9mm,-4mm>*{^{_2}}**@{},
   <-0.35mm,-0.35mm>*{};<-2.8mm,-4mm>*{^{_1}}**@{},
\end{xy}\Ea
=(-1)^{d}
\Ba{c}\begin{xy}
 <0mm,0.66mm>*{};<0mm,3mm>*{}**@{-},
 <0.39mm,-0.39mm>*{};<2.2mm,-2.2mm>*{}**@{-},
 <-0.35mm,-0.35mm>*{};<-2.2mm,-2.2mm>*{}**@{-},
 <0mm,0mm>*{\circ};<0mm,0mm>*{}**@{},
   <0.39mm,-0.39mm>*{};<2.9mm,-4mm>*{^{_1}}**@{},
   <-0.35mm,-0.35mm>*{};<-2.8mm,-4mm>*{^{_2}}**@{},
\end{xy}\Ea
\right\rangle
$$
by the ideal generated by the following elements
\Beq\label{R for LieB}
\cR:\left\{
\Ba{c}
\Ba{c}\resizebox{7mm}{!}{
\begin{xy}
 <0mm,0mm>*{\circ};<0mm,0mm>*{}**@{},
 <0mm,-0.49mm>*{};<0mm,-3.0mm>*{}**@{-},
 <0.49mm,0.49mm>*{};<1.9mm,1.9mm>*{}**@{-},
 <-0.5mm,0.5mm>*{};<-1.9mm,1.9mm>*{}**@{-},
 <-2.3mm,2.3mm>*{\circ};<-2.3mm,2.3mm>*{}**@{},
 <-1.8mm,2.8mm>*{};<0mm,4.9mm>*{}**@{-},
 <-2.8mm,2.9mm>*{};<-4.6mm,4.9mm>*{}**@{-},
   <0.49mm,0.49mm>*{};<2.7mm,2.3mm>*{^3}**@{},
   <-1.8mm,2.8mm>*{};<0.4mm,5.3mm>*{^2}**@{},
   <-2.8mm,2.9mm>*{};<-5.1mm,5.3mm>*{^1}**@{},
 \end{xy}}\Ea
 +
\Ba{c}\resizebox{7mm}{!}{\begin{xy}
 <0mm,0mm>*{\circ};<0mm,0mm>*{}**@{},
 <0mm,-0.49mm>*{};<0mm,-3.0mm>*{}**@{-},
 <0.49mm,0.49mm>*{};<1.9mm,1.9mm>*{}**@{-},
 <-0.5mm,0.5mm>*{};<-1.9mm,1.9mm>*{}**@{-},
 <-2.3mm,2.3mm>*{\circ};<-2.3mm,2.3mm>*{}**@{},
 <-1.8mm,2.8mm>*{};<0mm,4.9mm>*{}**@{-},
 <-2.8mm,2.9mm>*{};<-4.6mm,4.9mm>*{}**@{-},
   <0.49mm,0.49mm>*{};<2.7mm,2.3mm>*{^2}**@{},
   <-1.8mm,2.8mm>*{};<0.4mm,5.3mm>*{^1}**@{},
   <-2.8mm,2.9mm>*{};<-5.1mm,5.3mm>*{^3}**@{},
 \end{xy}}\Ea
 +
\Ba{c}\resizebox{7mm}{!}{\begin{xy}
 <0mm,0mm>*{\circ};<0mm,0mm>*{}**@{},
 <0mm,-0.49mm>*{};<0mm,-3.0mm>*{}**@{-},
 <0.49mm,0.49mm>*{};<1.9mm,1.9mm>*{}**@{-},
 <-0.5mm,0.5mm>*{};<-1.9mm,1.9mm>*{}**@{-},
 <-2.3mm,2.3mm>*{\circ};<-2.3mm,2.3mm>*{}**@{},
 <-1.8mm,2.8mm>*{};<0mm,4.9mm>*{}**@{-},
 <-2.8mm,2.9mm>*{};<-4.6mm,4.9mm>*{}**@{-},
   <0.49mm,0.49mm>*{};<2.7mm,2.3mm>*{^1}**@{},
   <-1.8mm,2.8mm>*{};<0.4mm,5.3mm>*{^3}**@{},
   <-2.8mm,2.9mm>*{};<-5.1mm,5.3mm>*{^2}**@{},
 \end{xy}}\Ea
 \ \ , \ \
\Ba{c}\resizebox{8.4mm}{!}{ \begin{xy}
 <0mm,0mm>*{\circ};<0mm,0mm>*{}**@{},
 <0mm,0.69mm>*{};<0mm,3.0mm>*{}**@{-},
 <0.39mm,-0.39mm>*{};<2.4mm,-2.4mm>*{}**@{-},
 <-0.35mm,-0.35mm>*{};<-1.9mm,-1.9mm>*{}**@{-},
 <-2.4mm,-2.4mm>*{\circ};<-2.4mm,-2.4mm>*{}**@{},
 <-2.0mm,-2.8mm>*{};<0mm,-4.9mm>*{}**@{-},
 <-2.8mm,-2.9mm>*{};<-4.7mm,-4.9mm>*{}**@{-},
    <0.39mm,-0.39mm>*{};<3.3mm,-4.0mm>*{^3}**@{},
    <-2.0mm,-2.8mm>*{};<0.5mm,-6.7mm>*{^2}**@{},
    <-2.8mm,-2.9mm>*{};<-5.2mm,-6.7mm>*{^1}**@{},
 \end{xy}}\Ea
 +
\Ba{c}\resizebox{8.4mm}{!}{ \begin{xy}
 <0mm,0mm>*{\circ};<0mm,0mm>*{}**@{},
 <0mm,0.69mm>*{};<0mm,3.0mm>*{}**@{-},
 <0.39mm,-0.39mm>*{};<2.4mm,-2.4mm>*{}**@{-},
 <-0.35mm,-0.35mm>*{};<-1.9mm,-1.9mm>*{}**@{-},
 <-2.4mm,-2.4mm>*{\circ};<-2.4mm,-2.4mm>*{}**@{},
 <-2.0mm,-2.8mm>*{};<0mm,-4.9mm>*{}**@{-},
 <-2.8mm,-2.9mm>*{};<-4.7mm,-4.9mm>*{}**@{-},
    <0.39mm,-0.39mm>*{};<3.3mm,-4.0mm>*{^2}**@{},
    <-2.0mm,-2.8mm>*{};<0.5mm,-6.7mm>*{^1}**@{},
    <-2.8mm,-2.9mm>*{};<-5.2mm,-6.7mm>*{^3}**@{},
 \end{xy}}\Ea
 +
\Ba{c}\resizebox{8.4mm}{!}{ \begin{xy}
 <0mm,0mm>*{\circ};<0mm,0mm>*{}**@{},
 <0mm,0.69mm>*{};<0mm,3.0mm>*{}**@{-},
 <0.39mm,-0.39mm>*{};<2.4mm,-2.4mm>*{}**@{-},
 <-0.35mm,-0.35mm>*{};<-1.9mm,-1.9mm>*{}**@{-},
 <-2.4mm,-2.4mm>*{\circ};<-2.4mm,-2.4mm>*{}**@{},
 <-2.0mm,-2.8mm>*{};<0mm,-4.9mm>*{}**@{-},
 <-2.8mm,-2.9mm>*{};<-4.7mm,-4.9mm>*{}**@{-},
    <0.39mm,-0.39mm>*{};<3.3mm,-4.0mm>*{^1}**@{},
    <-2.0mm,-2.8mm>*{};<0.5mm,-6.7mm>*{^3}**@{},
    <-2.8mm,-2.9mm>*{};<-5.2mm,-6.7mm>*{^2}**@{},
 \end{xy}}\Ea
 \\
 \Ba{c}\resizebox{5mm}{!}{\begin{xy}
 <0mm,2.47mm>*{};<0mm,0.12mm>*{}**@{-},
 <0.5mm,3.5mm>*{};<2.2mm,5.2mm>*{}**@{-},
 <-0.48mm,3.48mm>*{};<-2.2mm,5.2mm>*{}**@{-},
 <0mm,3mm>*{\circ};<0mm,3mm>*{}**@{},
  <0mm,-0.8mm>*{\circ};<0mm,-0.8mm>*{}**@{},
<-0.39mm,-1.2mm>*{};<-2.2mm,-3.5mm>*{}**@{-},
 <0.39mm,-1.2mm>*{};<2.2mm,-3.5mm>*{}**@{-},
     <0.5mm,3.5mm>*{};<2.8mm,5.7mm>*{^2}**@{},
     <-0.48mm,3.48mm>*{};<-2.8mm,5.7mm>*{^1}**@{},
   <0mm,-0.8mm>*{};<-2.7mm,-5.2mm>*{^1}**@{},
   <0mm,-0.8mm>*{};<2.7mm,-5.2mm>*{^2}**@{},
\end{xy}}\Ea
  -
\Ba{c}\resizebox{7mm}{!}{\begin{xy}
 <0mm,-1.3mm>*{};<0mm,-3.5mm>*{}**@{-},
 <0.38mm,-0.2mm>*{};<2.0mm,2.0mm>*{}**@{-},
 <-0.38mm,-0.2mm>*{};<-2.2mm,2.2mm>*{}**@{-},
<0mm,-0.8mm>*{\circ};<0mm,0.8mm>*{}**@{},
 <2.4mm,2.4mm>*{\circ};<2.4mm,2.4mm>*{}**@{},
 <2.77mm,2.0mm>*{};<4.4mm,-0.8mm>*{}**@{-},
 <2.4mm,3mm>*{};<2.4mm,5.2mm>*{}**@{-},
     <0mm,-1.3mm>*{};<0mm,-5.3mm>*{^1}**@{},
     <2.5mm,2.3mm>*{};<5.1mm,-2.6mm>*{^2}**@{},
    <2.4mm,2.5mm>*{};<2.4mm,5.7mm>*{^2}**@{},
    <-0.38mm,-0.2mm>*{};<-2.8mm,2.5mm>*{^1}**@{},
    \end{xy}}\Ea
  - (-1)^{d}
\Ba{c}\resizebox{7mm}{!}{\begin{xy}
 <0mm,-1.3mm>*{};<0mm,-3.5mm>*{}**@{-},
 <0.38mm,-0.2mm>*{};<2.0mm,2.0mm>*{}**@{-},
 <-0.38mm,-0.2mm>*{};<-2.2mm,2.2mm>*{}**@{-},
<0mm,-0.8mm>*{\circ};<0mm,0.8mm>*{}**@{},
 <2.4mm,2.4mm>*{\circ};<2.4mm,2.4mm>*{}**@{},
 <2.77mm,2.0mm>*{};<4.4mm,-0.8mm>*{}**@{-},
 <2.4mm,3mm>*{};<2.4mm,5.2mm>*{}**@{-},
     <0mm,-1.3mm>*{};<0mm,-5.3mm>*{^2}**@{},
     <2.5mm,2.3mm>*{};<5.1mm,-2.6mm>*{^1}**@{},
    <2.4mm,2.5mm>*{};<2.4mm,5.7mm>*{^2}**@{},
    <-0.38mm,-0.2mm>*{};<-2.8mm,2.5mm>*{^1}**@{},
    \end{xy}}\Ea
  - (-1)^{d+c}
\Ba{c}\resizebox{7mm}{!}{\begin{xy}
 <0mm,-1.3mm>*{};<0mm,-3.5mm>*{}**@{-},
 <0.38mm,-0.2mm>*{};<2.0mm,2.0mm>*{}**@{-},
 <-0.38mm,-0.2mm>*{};<-2.2mm,2.2mm>*{}**@{-},
<0mm,-0.8mm>*{\circ};<0mm,0.8mm>*{}**@{},
 <2.4mm,2.4mm>*{\circ};<2.4mm,2.4mm>*{}**@{},
 <2.77mm,2.0mm>*{};<4.4mm,-0.8mm>*{}**@{-},
 <2.4mm,3mm>*{};<2.4mm,5.2mm>*{}**@{-},
     <0mm,-1.3mm>*{};<0mm,-5.3mm>*{^2}**@{},
     <2.5mm,2.3mm>*{};<5.1mm,-2.6mm>*{^1}**@{},
    <2.4mm,2.5mm>*{};<2.4mm,5.7mm>*{^1}**@{},
    <-0.38mm,-0.2mm>*{};<-2.8mm,2.5mm>*{^2}**@{},
    \end{xy}}\Ea
 - (-1)^{c}
\Ba{c}\resizebox{7mm}{!}{\begin{xy}
 <0mm,-1.3mm>*{};<0mm,-3.5mm>*{}**@{-},
 <0.38mm,-0.2mm>*{};<2.0mm,2.0mm>*{}**@{-},
 <-0.38mm,-0.2mm>*{};<-2.2mm,2.2mm>*{}**@{-},
<0mm,-0.8mm>*{\circ};<0mm,0.8mm>*{}**@{},
 <2.4mm,2.4mm>*{\circ};<2.4mm,2.4mm>*{}**@{},
 <2.77mm,2.0mm>*{};<4.4mm,-0.8mm>*{}**@{-},
 <2.4mm,3mm>*{};<2.4mm,5.2mm>*{}**@{-},
     <0mm,-1.3mm>*{};<0mm,-5.3mm>*{^1}**@{},
     <2.5mm,2.3mm>*{};<5.1mm,-2.6mm>*{^2}**@{},
    <2.4mm,2.5mm>*{};<2.4mm,5.7mm>*{^1}**@{},
    <-0.38mm,-0.2mm>*{};<-2.8mm,2.5mm>*{^2}**@{},
    \end{xy}}\Ea
    \Ea
\right.
\Eeq

Similarly,  $\LoBcd$ (with $c+d\in 2\Z$ by default) is a quadratic properad
$
\cF ree\langle E\rangle/\langle\cR_\diamond\rangle
$
generated by the same $\bS$-bimodule $E$ modulo the relations
$$
\cR_\diamond:= \cR \ \bigsqcup
\Ba{c}\resizebox{4mm}{!}
{\xy
 (0,0)*{\circ}="a",
(0,6)*{\circ}="b",
(3,3)*{}="c",
(-3,3)*{}="d",
 (0,9)*{}="b'",
(0,-3)*{}="a'",
\ar@{-} "a";"c" <0pt>
\ar @{-} "a";"d" <0pt>
\ar @{-} "a";"a'" <0pt>
\ar @{-} "b";"c" <0pt>
\ar @{-} "b";"d" <0pt>
\ar @{-} "b";"b'" <0pt>
\endxy}
\Ea
$$

It is clear from the association
$
\vartriangle \leftrightarrow
 \begin{xy}
 <0mm,-0.55mm>*{};<0mm,-2.5mm>*{}**@{-},
 <0.5mm,0.5mm>*{};<2.2mm,2.2mm>*{}**@{-},
 <-0.48mm,0.48mm>*{};<-2.2mm,2.2mm>*{}**@{-},
 <0mm,0mm>*{\circ};<0mm,0mm>*{}**@{},
 \end{xy}$,
 $
[\ , \ ] \leftrightarrow
 \begin{xy}
 <0mm,0.66mm>*{};<0mm,3mm>*{}**@{-},
 <0.39mm,-0.39mm>*{};<2.2mm,-2.2mm>*{}**@{-},
 <-0.35mm,-0.35mm>*{};<-2.2mm,-2.2mm>*{}**@{-},
 <0mm,0mm>*{\circ};<0mm,0mm>*{}**@{},
 \end{xy}
$
that there is a one-to-one correspondence between representations of $\LBcd$ (resp.,
$\LoBcd$) in a finite dimensional space $V$ and (resp., involutive) Lie $(c+d-2)$-bialgebra
structures in $V[c-1]$.

\sip

The minimal resolution $\HoLBcd$ of the properad $\LBcd$ was constructed in \cite{Ko,MaVo} for $d+c\in 2\Z$ and in
 \cite{Me1,Me2} for $d+c\in 2\Z+1$. It
is generated by the following (skew)symmetric corollas of degree $1 +c(1-m)+d(1-n)$
\Beq\label{2: symmetries of HoLiebcd corollas}
\Ba{c}\resizebox{17mm}{!}{\begin{xy}
 <0mm,0mm>*{\circ};<0mm,0mm>*{}**@{},
 <-0.6mm,0.44mm>*{};<-8mm,5mm>*{}**@{-},
 <-0.4mm,0.7mm>*{};<-4.5mm,5mm>*{}**@{-},
 <0mm,0mm>*{};<1mm,5mm>*{\ldots}**@{},
 <0.4mm,0.7mm>*{};<4.5mm,5mm>*{}**@{-},
 <0.6mm,0.44mm>*{};<8mm,5mm>*{}**@{-},
   <0mm,0mm>*{};<-10.5mm,5.9mm>*{^{\sigma(1)}}**@{},
   <0mm,0mm>*{};<-4mm,5.9mm>*{^{\sigma(2)}}**@{},
   <0mm,0mm>*{};<10.0mm,5.9mm>*{^{\sigma(m)}}**@{},
 <-0.6mm,-0.44mm>*{};<-8mm,-5mm>*{}**@{-},
 <-0.4mm,-0.7mm>*{};<-4.5mm,-5mm>*{}**@{-},
 <0mm,0mm>*{};<1mm,-5mm>*{\ldots}**@{},
 <0.4mm,-0.7mm>*{};<4.5mm,-5mm>*{}**@{-},
 <0.6mm,-0.44mm>*{};<8mm,-5mm>*{}**@{-},
   <0mm,0mm>*{};<-10.5mm,-6.9mm>*{^{\tau(1)}}**@{},
   <0mm,0mm>*{};<-4mm,-6.9mm>*{^{\tau(2)}}**@{},
   <0mm,0mm>*{};<10.0mm,-6.9mm>*{^{\tau(n)}}**@{},
 \end{xy}}\Ea
=(-1)^{c|\sigma|+d|\tau|}
\Ba{c}\resizebox{14mm}{!}{\begin{xy}
 <0mm,0mm>*{\circ};<0mm,0mm>*{}**@{},
 <-0.6mm,0.44mm>*{};<-8mm,5mm>*{}**@{-},
 <-0.4mm,0.7mm>*{};<-4.5mm,5mm>*{}**@{-},
 <0mm,0mm>*{};<-1mm,5mm>*{\ldots}**@{},
 <0.4mm,0.7mm>*{};<4.5mm,5mm>*{}**@{-},
 <0.6mm,0.44mm>*{};<8mm,5mm>*{}**@{-},
   <0mm,0mm>*{};<-8.5mm,5.5mm>*{^1}**@{},
   <0mm,0mm>*{};<-5mm,5.5mm>*{^2}**@{},
   <0mm,0mm>*{};<4.5mm,5.5mm>*{^{m\hspace{-0.5mm}-\hspace{-0.5mm}1}}**@{},
   <0mm,0mm>*{};<9.0mm,5.5mm>*{^m}**@{},
 <-0.6mm,-0.44mm>*{};<-8mm,-5mm>*{}**@{-},
 <-0.4mm,-0.7mm>*{};<-4.5mm,-5mm>*{}**@{-},
 <0mm,0mm>*{};<-1mm,-5mm>*{\ldots}**@{},
 <0.4mm,-0.7mm>*{};<4.5mm,-5mm>*{}**@{-},
 <0.6mm,-0.44mm>*{};<8mm,-5mm>*{}**@{-},
   <0mm,0mm>*{};<-8.5mm,-6.9mm>*{^1}**@{},
   <0mm,0mm>*{};<-5mm,-6.9mm>*{^2}**@{},
   <0mm,0mm>*{};<4.5mm,-6.9mm>*{^{n\hspace{-0.5mm}-\hspace{-0.5mm}1}}**@{},
   <0mm,0mm>*{};<9.0mm,-6.9mm>*{^n}**@{},
 \end{xy}}\Ea \ \ \forall \sigma\in \bS_m, \forall\tau\in \bS_n
\Eeq
and has the differential
given on the generators by
\Beq\label{LBk_infty}
\delta
\Ba{c}\resizebox{14mm}{!}{\begin{xy}
 <0mm,0mm>*{\circ};<0mm,0mm>*{}**@{},
 <-0.6mm,0.44mm>*{};<-8mm,5mm>*{}**@{-},
 <-0.4mm,0.7mm>*{};<-4.5mm,5mm>*{}**@{-},
 <0mm,0mm>*{};<-1mm,5mm>*{\ldots}**@{},
 <0.4mm,0.7mm>*{};<4.5mm,5mm>*{}**@{-},
 <0.6mm,0.44mm>*{};<8mm,5mm>*{}**@{-},
   <0mm,0mm>*{};<-8.5mm,5.5mm>*{^1}**@{},
   <0mm,0mm>*{};<-5mm,5.5mm>*{^2}**@{},
   <0mm,0mm>*{};<4.5mm,5.5mm>*{^{m\hspace{-0.5mm}-\hspace{-0.5mm}1}}**@{},
   <0mm,0mm>*{};<9.0mm,5.5mm>*{^m}**@{},
 <-0.6mm,-0.44mm>*{};<-8mm,-5mm>*{}**@{-},
 <-0.4mm,-0.7mm>*{};<-4.5mm,-5mm>*{}**@{-},
 <0mm,0mm>*{};<-1mm,-5mm>*{\ldots}**@{},
 <0.4mm,-0.7mm>*{};<4.5mm,-5mm>*{}**@{-},
 <0.6mm,-0.44mm>*{};<8mm,-5mm>*{}**@{-},
   <0mm,0mm>*{};<-8.5mm,-6.9mm>*{^1}**@{},
   <0mm,0mm>*{};<-5mm,-6.9mm>*{^2}**@{},
   <0mm,0mm>*{};<4.5mm,-6.9mm>*{^{n\hspace{-0.5mm}-\hspace{-0.5mm}1}}**@{},
   <0mm,0mm>*{};<9.0mm,-6.9mm>*{^n}**@{},
 \end{xy}}\Ea
\ \ = \ \
 \sum_{[1,\ldots,m]=I_1\sqcup I_2\atop
 {|I_1|\geq 0, |I_2|\geq 1}}
 \sum_{[1,\ldots,n]=J_1\sqcup J_2\atop
 {|J_1|\geq 1, |J_2|\geq 1}
}\hspace{0mm}
\pm
\Ba{c}\resizebox{22mm}{!}{ \begin{xy}
 <0mm,0mm>*{\circ};<0mm,0mm>*{}**@{},
 <-0.6mm,0.44mm>*{};<-8mm,5mm>*{}**@{-},
 <-0.4mm,0.7mm>*{};<-4.5mm,5mm>*{}**@{-},
 <0mm,0mm>*{};<0mm,5mm>*{\ldots}**@{},
 <0.4mm,0.7mm>*{};<4.5mm,5mm>*{}**@{-},
 <0.6mm,0.44mm>*{};<12.4mm,4.8mm>*{}**@{-},
     <0mm,0mm>*{};<-2mm,7mm>*{\overbrace{\ \ \ \ \ \ \ \ \ \ \ \ }}**@{},
     <0mm,0mm>*{};<-2mm,9mm>*{^{I_1}}**@{},
 <-0.6mm,-0.44mm>*{};<-8mm,-5mm>*{}**@{-},
 <-0.4mm,-0.7mm>*{};<-4.5mm,-5mm>*{}**@{-},
 <0mm,0mm>*{};<-1mm,-5mm>*{\ldots}**@{},
 <0.4mm,-0.7mm>*{};<4.5mm,-5mm>*{}**@{-},
 <0.6mm,-0.44mm>*{};<8mm,-5mm>*{}**@{-},
      <0mm,0mm>*{};<0mm,-7mm>*{\underbrace{\ \ \ \ \ \ \ \ \ \ \ \ \ \ \
      }}**@{},
      <0mm,0mm>*{};<0mm,-10.6mm>*{_{J_1}}**@{},
 <13mm,5mm>*{};<13mm,5mm>*{\circ}**@{},
 <12.6mm,5.44mm>*{};<5mm,10mm>*{}**@{-},
 <12.6mm,5.7mm>*{};<8.5mm,10mm>*{}**@{-},
 <13mm,5mm>*{};<13mm,10mm>*{\ldots}**@{},
 <13.4mm,5.7mm>*{};<16.5mm,10mm>*{}**@{-},
 <13.6mm,5.44mm>*{};<20mm,10mm>*{}**@{-},
      <13mm,5mm>*{};<13mm,12mm>*{\overbrace{\ \ \ \ \ \ \ \ \ \ \ \ \ \ }}**@{},
      <13mm,5mm>*{};<13mm,14mm>*{^{I_2}}**@{},
 <12.4mm,4.3mm>*{};<8mm,0mm>*{}**@{-},
 <12.6mm,4.3mm>*{};<12mm,0mm>*{\ldots}**@{},
 <13.4mm,4.5mm>*{};<16.5mm,0mm>*{}**@{-},
 <13.6mm,4.8mm>*{};<20mm,0mm>*{}**@{-},
     <13mm,5mm>*{};<14.3mm,-2mm>*{\underbrace{\ \ \ \ \ \ \ \ \ \ \ }}**@{},
     <13mm,5mm>*{};<14.3mm,-4.5mm>*{_{J_2}}**@{},
 \end{xy}}\Ea
\Eeq
where the signs on the r.h.s\ are uniquely fixed for $c+d\in 2\Z$ by the fact that they all equal to $+1$ if $ c$ and $d$ are even integers, and for $c+d\in 2\Z+1$ the signs are given explicitly in
\cite{Me1}.

\sip

The minimal resolution $\HoLoBcd$ of the properad $\LoBcd$ was constructed in \cite{CMW}. It is a free properad generated
by the following (skew)symmetric corollas of degree $1+c(1-m-a)+d(1-n-a)$
\Beq\label{equ:LoBgenerators}
\Ba{c}\resizebox{16mm}{!}{\xy
(-9,-6)*{};
(0,0)*+{a}*\cir{}
**\dir{-};
(-5,-6)*{};
(0,0)*+{a}*\cir{}
**\dir{-};
(9,-6)*{};
(0,0)*+{a}*\cir{}
**\dir{-};
(5,-6)*{};
(0,0)*+{a}*\cir{}
**\dir{-};
(0,-6)*{\ldots};
(-10,-8)*{_1};
(-6,-8)*{_2};
(10,-8)*{_n};
(-9,6)*{};
(0,0)*+{a}*\cir{}
**\dir{-};
(-5,6)*{};
(0,0)*+{a}*\cir{}
**\dir{-};
(9,6)*{};
(0,0)*+{a}*\cir{}
**\dir{-};
(5,6)*{};
(0,0)*+{a}*\cir{}
**\dir{-};
(0,6)*{\ldots};
(-10,8)*{_1};
(-6,8)*{_2};
(10,8)*{_m};
\endxy}\Ea
=(-1)^{(d+1)(\sigma+\tau)}
\Ba{c}\resizebox{20mm}{!}{\xy
(-9,-6)*{};
(0,0)*+{a}*\cir{}
**\dir{-};
(-5,-6)*{};
(0,0)*+{a}*\cir{}
**\dir{-};
(9,-6)*{};
(0,0)*+{a}*\cir{}
**\dir{-};
(5,-6)*{};
(0,0)*+{a}*\cir{}
**\dir{-};
(0,-6)*{\ldots};
(-12,-8)*{_{\tau(1)}};
(-6,-8)*{_{\tau(2)}};
(12,-8)*{_{\tau(n)}};
(-9,6)*{};
(0,0)*+{a}*\cir{}
**\dir{-};
(-5,6)*{};
(0,0)*+{a}*\cir{}
**\dir{-};
(9,6)*{};
(0,0)*+{a}*\cir{}
**\dir{-};
(5,6)*{};
(0,0)*+{a}*\cir{}
**\dir{-};
(0,6)*{\ldots};
(-12,8)*{_{\sigma(1)}};
(-6,8)*{_{\sigma(2)}};
(12,8)*{_{\sigma(m)}};
\endxy}\Ea\ \ \ \forall \sigma\in \bS_m, \forall \tau\in \bS_n,
\Eeq
where $m+n+ a\geq 3$, $m\geq 1$, $n\geq 1$, $a\geq 0$. The differential in
$\HoLoBd$ is given on the generators by
\Beq\label{2: d on Lie inv infty}
\delta
\Ba{c}\resizebox{16mm}{!}{\xy
(-9,-6)*{};
(0,0)*+{a}*\cir{}
**\dir{-};
(-5,-6)*{};
(0,0)*+{a}*\cir{}
**\dir{-};
(9,-6)*{};
(0,0)*+{a}*\cir{}
**\dir{-};
(5,-6)*{};
(0,0)*+{a}*\cir{}
**\dir{-};
(0,-6)*{\ldots};
(-10,-8)*{_1};
(-6,-8)*{_2};
(10,-8)*{_n};
(-9,6)*{};
(0,0)*+{a}*\cir{}
**\dir{-};
(-5,6)*{};
(0,0)*+{a}*\cir{}
**\dir{-};
(9,6)*{};
(0,0)*+{a}*\cir{}
**\dir{-};
(5,6)*{};
(0,0)*+{a}*\cir{}
**\dir{-};
(0,6)*{\ldots};
(-10,8)*{_1};
(-6,8)*{_2};
(10,8)*{_m};
\endxy}\Ea
=
\sum_{l\geq 1}\sum_{a=b+c+l-1}\sum_{[m]=I_1\sqcup I_2\atop
[n]=J_1\sqcup J_2} \pm
\Ba{c}
%
%
\Ba{c}\resizebox{21mm}{!}{\xy
(0,0)*+{b}*\cir{}="b",
(10,10)*+{c}*\cir{}="c",
%
(-9,6)*{}="1",
(-7,6)*{}="2",
(-2,6)*{}="3",
(-3.5,5)*{...},
(-4,-6)*{}="-1",
(-2,-6)*{}="-2",
(4,-6)*{}="-3",
(1,-5)*{...},
(0,-8)*{\underbrace{\ \ \ \ \ \ \ \ }},
(0,-11)*{_{J_1}},
(-6,8)*{\overbrace{ \ \ \ \ \ \ }},
(-6,11)*{_{I_1}},
(6,16)*{}="1'",
(8,16)*{}="2'",
(14,16)*{}="3'",
(11,15)*{...},
(11,6)*{}="-1'",
(16,6)*{}="-2'",
(18,6)*{}="-3'",
(13.5,6)*{...},
(15,4)*{\underbrace{\ \ \ \ \ \ \ }},
(15,1)*{_{J_2}},
(10,18)*{\overbrace{ \ \ \ \ \ \ \ \ }},
(10,21)*{_{I_2}},
%
(0,2)*-{};(8.0,10.0)*-{}
**\crv{(0,10)};
(0.5,1.8)*-{};(8.5,9.0)*-{}
**\crv{(0.4,7)};
(1.5,0.5)*-{};(9.1,8.5)*-{}
**\crv{(5,1)};
(1.7,0.0)*-{};(9.5,8.6)*-{}
**\crv{(6,-1)};
(5,5)*+{...};
\ar @{-} "b";"1" <0pt>
\ar @{-} "b";"2" <0pt>
\ar @{-} "b";"3" <0pt>
\ar @{-} "b";"-1" <0pt>
\ar @{-} "b";"-2" <0pt>
\ar @{-} "b";"-3" <0pt>
\ar @{-} "c";"1'" <0pt>
\ar @{-} "c";"2'" <0pt>
\ar @{-} "c";"3'" <0pt>
\ar @{-} "c";"-1'" <0pt>
\ar @{-} "c";"-2'" <0pt>
\ar @{-} "c";"-3'" <0pt>
\endxy}\Ea
\Ea
\Eeq
where the summation parameter $l$ counts the number of internal edges connecting the two vertices
on the r.h.s., and the signs are  fixed by the fact that they all equal to $+1$ for $c$ and $d$
odd integers.

\sip

Our purpose in this paper is to relate deformation complexes of all the properads
considered above to various graph complexes whose cohomology is partially computed, and whose relations with the Grothendieck-Teuchm\"uller Lie algebra are well-understood.

\subsection{Complete variants}\label{sec:completions}
Note that the defining relations for the properads $\LBcd$ and $\LoBcd$  do not mix composition diagrams of different loop orders. It follow that the mentioned properads are all graded by the loop order (here also called genus) of composition diagrams.
In particular, fixing the arity, the operations are finite linear combinations (not series) of composites of generators.
For some applications, including in particular the integration of derivations to automorphisms, it is more convenient to consider the completed versions by the genus grading $\wLBcd$ and $\wLoBcd$. Concretely, the operations of fixed arity in the complete versions of our properads are given by infinite series (instead of just linear combinations) of composites of generators.

\sip

Similarly, the resolutions $\HoLBcd$ and $\HoLoBcd$  receive a grading by the loop order (or genus), and we may also consider the completed versions (with respect to this grading) $\wHoLBcd$ and $\wHoLoBcd$.
Here it should be noted that with respect to the genus grading the generator \eqref{equ:LoBgenerators} must be considered as living in degree $a$ to make this  grading consistent with the definition of differential \eqref{2: d on Lie inv infty}.

\subsection{Directed graph complexes} A {\em graph}\, $\Ga$ is a 1-dimensional $CW$ complex whose 0-cells are called {\em vertices}\, and 1-cells are called {\em edges}. The set of vertices of $\Ga$ is denoted by $V(\Ga)$ and the set of edges by $E(\Ga)$. A graph $\Ga$ is called {\em directed}\, if its edge $e\in E(\Ga)$ comes equipped with an orientation or, plainly speaking,
with a choice of a direction.

\sip

Let $G_{n,l}$ be a set of directed graphs $\Ga$ with $n$ vertices and $l$  edges such that
some bijections $V(\Ga)\rar [n]$ and $E(\Ga)\rar [l]$ are fixed, i.e.\ every edges and every vertex of $\Ga$ has a fixed numerical label. There is
a natural right action of the group $\bS_n \times  \bS_l$ on the set $G_{n,l}$ with $\bS_n$ acting by relabeling the vertices and  $\bS_l$ by relabeling the
edges. 

\sip

For each fixed integer $d$, a collection of $\bS_n$-modules,
$$
\caD\cG ra_{d}=\left\{\caD\cG ra_d(n):= \prod_{l\geq 0} \K \langle G_{n,l}\rangle \ot_{ \bS_l}  \sgn_l^{\ot |d-1|} [l(d-1)]   \right\}_{n\geq 1}
$$
 is an operad with respect to the following operadic composition,
$$
\Ba{rccc}
\circ_i: &  \caD\cG ra_d(n) \times \caD\cG ra_d(m) &\lon & \caD\cG ra_d(m+n-1),  \ \ \forall\ i\in [n]\\
         &       (\Ga_1, \Ga_2) &\lon &      \Ga_1\circ_i \Ga_2,
\Ea
$$
where  $\Ga_1\circ_i \Ga_2$ is defined by substituting the graph $\Ga_2$ into the $i$-labeled vertex $v_i$ of $\Ga_1$ and taking a sum over  re-attachments of dangling edges (attached before to $v_i$) to vertices of $\Ga_2$
in all possible ways.

\sip

For any operad $\cP=\{\cP(n)\}_{n\geq 1}$  in the category of graded vector spaces,
the linear the map
$$
\Ba{rccc}
[\ ,\ ]:&  \sP \ot \sP & \lon & \sP\\
& (a\in \cP(n), b\in \cP(m)) & \lon &
[a, b]:= \sum_{i=1}^n a\circ_i b - (-1)^{|a||b|}\sum_{i=1}^m b\circ_i a\ \in \cP(m+n-1)
\Ea
$$
makes a graded vector space
$
\sP:= \prod_{n\geq 1}\cP(n)$
into a Lie algebra \cite{KM}; moreover, these brackets induce a Lie algebra structure on the subspace
of invariants
$
\sP^\bS:=  \prod_{n\geq 1}\cP(n)^{\bS_n}$. In particular,
the graded vector space
$$
\mathsf{dfGC}_{d}:= \prod_{n\geq 1} \cG ra_{d}(n)^{\bS_n}[d(1-n)]
$$
is a Lie algebra with respect to the above Lie brackets, and as such it can be identified
with the deformation complex $\Def(\caL ie_d\stackrel{0}{\rar} \cG ra_{d})$ of a zero morphism. Hence non-trivial Maurer-Cartan elements of $(\mathsf{fGC}_{d}, [\ ,\ ])$ give us non-trivial morphisms of operads
$$
f:\caL ie_d {\lon} \caD\cG ra_{d}.
$$
 One such non-trivial morphism $f$ is given explicitly on the generator of $\caL ie_{d}$ by \cite{Wi1}
\Beq\label{2:  map from Lie to dgra}
f \left(\Ba{c}\begin{xy}
 <0mm,0.66mm>*{};<0mm,3mm>*{}**@{-},
 <0.39mm,-0.39mm>*{};<2.2mm,-2.2mm>*{}**@{-},
 <-0.35mm,-0.35mm>*{};<-2.2mm,-2.2mm>*{}**@{-},
 <0mm,0mm>*{\circ};<0mm,0mm>*{}**@{},
   <0.39mm,-0.39mm>*{};<2.9mm,-4mm>*{^{_2}}**@{},
   <-0.35mm,-0.35mm>*{};<-2.8mm,-4mm>*{^{_1}}**@{},
\end{xy}\Ea\right)=
\Ba{c}\resizebox{6.3mm}{!}{\xy
(0,1)*+{_1}*\cir{}="b",
(8,1)*+{_2}*\cir{}="c",
\ar @{->} "b";"c" <0pt>
\endxy}
\Ea  - (-1)^d
\Ba{c}\resizebox{7mm}{!}{\xy
(0,1)*+{_2}*\cir{}="b",
(8,1)*+{_1}*\cir{}="c",
\ar @{->} "b";"c" <0pt>
\endxy}
\Ea=:\xy
 (0,0)*{\bullet}="a",
(5,0)*{\bu}="b",
\ar @{->} "a";"b" <0pt>
\endxy
\Eeq
Note that elements of $\mathsf{dfGC}_{d}$ can be identified with graphs from $\caD\cG ra_d$ whose vertices' labels are symmetrized (for $d$ even) or skew-symmetrized (for $d$ odd) so that in pictures we can forget about labels of vertices  and denote them by unlabelled black bullets as in the formula above. Note also that graphs from  $\mathsf{dfGC}_{d}$ come equipped with an orientation, $or$, which is a choice of ordering of edges (for $d$ even) or a choice of ordering of vertices (for $d$ odd) up to an even permutation on both cases. Thus every graph $\Ga\in \mathsf{dfGC}_{d}$  has at most two different orientations, $pr$ and $or^{opp}$, and one has
the standard relation, $(\Ga, or)=-(\Ga, or^{opp})$; as usual, the data $(\Ga, or)$ is abbreviate by $\Ga$ (with some choice of orientation implicitly assumed).  Note that the homological degree of graph $\Ga$ from $\mathsf{dfGC}_{d}$ is given by
$
|\Ga|=d(\# V(\Ga) -1) + (1-d) \# E(\Ga).
$

\sip

The above morphism (\ref{2:  map from Lie to dgra}) makes
 $(\mathsf{dfGC}_{d}, [\ ,\ ])$ into a {\em differential}\, Lie algebra with the differential
 $$
 \delta:= [\xy
 (0,0)*{\bullet}="a",
(5,0)*{\bu}="b",
\ar @{->} "a";"b" <0pt>
\endxy ,\ ].
 $$
 This dg Lie algebra contains a  dg subalgebra $\mathsf{dGC}_{d}$ spanned by connected graphs
with at least bivalent vertices.
%
%
It was proven in \cite{Wi1} that
$$
H^\bu(\mathsf{dfGC}_{d})= \mathsf{dGC}_{d}
$$
so that there is no loss of generality of working with $\mathsf{dGC}_{d}$ instead of
$\mathsf{dfGC}_{d}$. Moreover, one has an isomorphism of Lie algebras \cite{Wi1},
$$
H^0(\mathsf{dGC}_{d})=\fg\fr\ft_1,
$$
where $\fg\fr\ft_1$ is the Lie algebra of the Grothendieck-teichm\"u ller group $GRT_1$ introduced by Drinfeld in the context of deformation quantization of Lie bialgebras. Nowadays, this group play an important role in many  other areas of mathematics (e.g.\ in the knot theory, in deformation quantization
of Poisson manifolds,  and in the classification theory of solutions of Kashiwara-Vergne problem, see \cite{Fu} for a review, and many references cited there).

\subsubsection{\bf Remark} Often one considers instead of $\caD\cG ra_d$ an operad $\cG ra_d=\{\cG ra_d(n)\}$ defined by
$$
\cG ra_d(n):= \prod_{l\geq 0} \K \langle G_{n,l}\rangle \ot_{ \bS_l \ltimes
 (\bS_2)^l}  \sgn_l^{|d|}\ot \sgn_2^{\ot l|d-1|} [l(d-1)]
$$
where the group $(\bS_2)^l$ acts on graphs from $G_{n,l}$ by flipping directions of the edges.
Then, arguing as above, one arrives to the  graph complex
$$
\mathsf{fGC}_d:=\Def(\caL ie_d\stackrel{f}{\rar} \cG ra_{d})
$$
of {\em undirected}\, graphs. It contains three important dg Lie subalgebras:
(i) $\mathsf{fcGC}_d \subset \mathsf{fGC}_d$ which is spanned by connected graphs, (ii) $\mathsf{GC}_d^2 \subset \mathsf{fcGC}_d$ which is spanned by graphs with at least bivalent vertices, and (iii) $\mathsf{GC}_d \subset \mathsf{GC}_d^2$ which is spanned by graphs with at least trivalent vertices.
It was shown in \cite{Ko1,Wi1} that the cohomology of these subalgebras (which determine completely the cohomology  of the full graph complex $\mathsf{fGC}_d$) are related to each other as follows
$$
H^\bu(\mathsf{fGC}_{d}) = H^\bu(\mathsf{GC}_{d}^2) = H^\bu(\mathsf{GC}_{d})\ \oplus\ \bigoplus_{j\geq 1\atop j\equiv 2d+1 \mod 4} \K[d-j],
$$
where the summand  $ \K[d-j]$ is generated by the loop-type graph with $j$ binary vertices.
It was proven in \cite{Wi1} that the complex of directed graphs and the complex of undirected
graphs have the same cohomology,
$$
H^\bu(\mathsf{dGC}_{d})=H^\bu(\mathsf{GC}_{d}^2).
$$
In the present context it is more suitable to work with the directed complex $\mathsf{dGC}_{d}$
rather than with  $\mathsf{GC}_{d}^2$.

\subsection{Oriented graph complexes} A graphs $\Ga$ from the operad $\caD\cG ra_d$
is called {\em oriented}\, if it contains no {\em wheels}, that is, directed paths of edges
forming a closed circle. The subspace
$\cG ra_d^{or}\subset \caD \cG ra_d$ spanned by oriented graphs is a suboperad. For example,
$$
\Ba{c}
 \xy
(0,-1.5)*{_{_2}},
(6.2,-1.5)*{_{_1}},
(3,6)*{^{^3}},
 (0,0)*{\bullet}="a",
(6,0)*{\bu}="b",
(3,5)*{\bu}="c",
\ar @{->} "a";"b" <0pt>
\ar @{->} "a";"c" <0pt>
\ar @{<-} "c";"b" <0pt>
\endxy\Ea
\in \cG ra_d^{or} \ \ \ \ \mbox{ but}\ \ \ \
\Ba{c}
 \xy
(0,-1.5)*{_{_2}},
(6.2,-1.5)*{_{_1}},
(3,6)*{^{^3}},
 (0,0)*{\bullet}="a",
(6,0)*{\bu}="b",
(3,5)*{\bu}="c",
\ar @{->} "a";"b" <0pt>
\ar @{<-} "a";"c" <0pt>
\ar @{->} "c";"b" <0pt>
\endxy\Ea
\not\in \cG ra_d^{or}(3).
$$
The morphism (\ref{2:  map from Lie to dgra}) factors through the inclusion $\cG ra_d\subset
\caD\cG ra_d$ so that one can consider a graph complex
$$
\mathsf{fGC}^{or}_d:=\Def\left(\caL ie_d \stackrel{f}{\rar} \cG ra_d^{or}\right)
$$
and its subcomplex $\GCor_d$ spanned by connected graphs with at least bivalent vertices
and with no bivalent vertices of the form  $\xy
 (0,0)*{}="a",
(4,0)*{\bu}="b",
(8,0)*{}="c",
\ar @{->} "a";"b" <0pt>
\ar @{->} "b";"c" <0pt>
\endxy$. This subcomplex determines the cohomology of the full graph complex,
$H^\bu(\mathsf{fGC}^{or}_d)=\odot^\bu (H^\bu(\GCor_d)$.
It was proven in \cite{Wi2} that
$$
H^\bu(\GCor_{d+1})=H^\bu(\mathsf{dGC}_d)=H^\bu(\mathsf{GC}_d^2).
$$
In particular, one has a remarkable isomorphism of Lie algebras,
$
H^0(\GCor_3)=\fg\fr\ft$. Moreover $H^i(\GCor_3)=0$ for $i\leq 2$ and $H^{-1}(\GCor_3)$
is a 1-dimensional space generated by the graph $
\Ba{c}\resizebox{4mm}{!}{   \xy
   \ar@/^0.6pc/(0,-5)*{\bullet};(0,5)*{\bullet}
   \ar@/^{-0.6pc}/(0,-5)*{\bullet};(0,5)*{\bullet}
 \endxy}\Ea
$.

\mip

Consider next a Lie algebra $(\GCor_3[[\hbar]], [\ ,\ ]$, where $\GCor_3[[\hbar]]$ is the topological vector space spanned by formal power series in a formal parameter $\hbar$ of homological degree $2$, and $[\ ,\ ]$ are the Lie brackets obtained from the standard ones in
$\GCor_d$ by the continuous extension. It was shown in \cite{CMW} that the formal power series
\Beq\label{2: Phi_hbar MC element}
\Phi_\hbar:= \sum_{k=1}^\infty \hbar^{k-1} \underbrace{
\Ba{c}\resizebox{6mm}{!}  {\xy
(0,0)*{...},
   \ar@/^1pc/(0,-5)*{\bullet};(0,5)*{\bullet}
   \ar@/^{-1pc}/(0,-5)*{\bullet};(0,5)*{\bullet}
   \ar@/^0.6pc/(0,-5)*{\bullet};(0,5)*{\bullet}
   \ar@/^{-0.6pc}/(0,-5)*{\bullet};(0,5)*{\bullet}
 \endxy}
 \Ea}_{k\ \mathrm{edges}}
 \Eeq
is a Maurer-Cartan element in the Lie algebra $(\mathsf{f}\sG\sC_3^{or}[[\hbar]], [\ ,\ ]$ and hence makes the latter into a {\em differential}\, Lie algebra with the differential
$$
\delta_\hbar=[\Phi_\hbar,\ ].
$$
It was proven in \cite{CMW} that
 $H^0( \GC_3^{or}[[\hbar]], \delta_\hbar)\simeq
H^0(\GC_3^{or},\delta)\simeq \grt_1$ as Lie algebras. Moreover, $H^i( \GC_3^{or}[[\hbar]], \delta_\hbar)=0$
for all $i\leq -2$ and $H^{-1}( \GC_3^{or}[[\hbar]], \delta_\hbar)$ is a 1-dimensional vector space  class generated by the formal power series
$
\sum_{k=2}^\infty (k-1)\hbar^{k-2}\underbrace{
\Ba{c}\resizebox{6mm}{!}  {\xy
(0,0)*{...},
   \ar@/^1pc/(0,-5)*{\bullet};(0,5)*{\bullet}
   \ar@/^{-1pc}/(0,-5)*{\bullet};(0,5)*{\bullet}
   \ar@/^0.6pc/(0,-5)*{\bullet};(0,5)*{\bullet}
   \ar@/^{-0.6pc}/(0,-5)*{\bullet};(0,5)*{\bullet}
 \endxy}
 \Ea}_{k\ \mathrm{edges}}
$

Sometimes we do not show in our pictures directions of edges of oriented graphs
by assuming tacitly that the flow goes from the bottom to the top (as in the case of properads).

\bip

{\Large
\section{\bf Deformation complexes of  properads\\ and directed graph complexes}
}

\mip

\subsection{Deformation complexes of properads}

 \label{2 sec:defcomplexes}
 For $\cC$ a coaugmented co(pr)operad, we will denote by $\Omega(\cC)$ its cobar construction.
Concretely, $\Omega(\cC)=\cF ree\langle\overline \cC[-1]\rangle$ as a graded (pr)operad where $\overline \cC$ the cokernel of the coaugmetation and $\cF ree\langle\dots\rangle$ denotes the free (pr)operad generated by an $\bS$-(bi)module.
We will often use complexes of derivations of (pr)operads and deformation complexes of (pr)operad maps.
For a map of properads $f: \Omega(\cC){\to} \cP$, we will denote by
\begin{equation}\label{equ:Defdefi}
\Def( \Omega(\cC)\stackrel{f}{\to} \cP )\cong \prod_{m,n} \Hom_{\bS_m\times \bS_n}(\cC(m,n), \cP(m,n))
\end{equation}
the associated convolution complex. It is natural structure of a dg Lie algebra \cite{MV} controlling deformations of the morphism $f$.

\sip

We will also consider the  Lie algebra  $\Der(\cP)$  of derivations of the properad $\cP$; in fact,
we will use a minor variation of the standard definition (given, e.g., in \cite{Ta}) defined as follows. Let $\cP^+$ be the free properad generated by $\cP$ and one other operation
$\begin{xy}
 <0mm,-0.55mm>*{};<0mm,-3mm>*{}**@{-},
 <0mm,0.5mm>*{};<0mm,3mm>*{}**@{-},
 <0mm,0mm>*{\bullet};<0mm,0mm>*{}**@{},
 \end{xy}$ of arity $(1,1)$ and of cohomological degree $+1$. On $\cP^+$ we define a differential $\delta^+$ by setting its value on the new generator by
 $$
 \delta^+ \begin{xy}
 <0mm,-0.55mm>*{};<0mm,-3mm>*{}**@{-},
 <0mm,0.5mm>*{};<0mm,3mm>*{}**@{-},
 <0mm,0mm>*{\bullet};<0mm,0mm>*{}**@{},
 \end{xy} := \begin{xy}
 <0mm,0mm>*{};<0mm,-3mm>*{}**@{-},
 <0mm,0mm>*{};<0mm,6mm>*{}**@{-},
 <0mm,0mm>*{\bullet};
 <0mm,3mm>*{\bullet};
 \end{xy}
 $$
 and on any other element $a\in \cP(m,n)$ (which we identify pictorially with the $(m,n)$-corolla
 whose vertex is decorated with $a$) by the formula
 $$
 \delta^+
 \begin{xy}
 <0mm,0mm>*{\bullet};<0mm,0mm>*{}**@{},
 <0mm,0mm>*{};<-8mm,5mm>*{}**@{-},
 <0mm,0mm>*{};<-4.5mm,5mm>*{}**@{-},
 <0mm,0mm>*{};<-1mm,5mm>*{\ldots}**@{},
 <0mm,0mm>*{};<4.5mm,5mm>*{}**@{-},
 <0mm,0mm>*{};<8mm,5mm>*{}**@{-},
   <0mm,0mm>*{};<-8.5mm,5.5mm>*{^1}**@{},
   <0mm,0mm>*{};<-5mm,5.5mm>*{^2}**@{},
   <0mm,0mm>*{};<4.5mm,5.5mm>*{^{m\hspace{-0.5mm}-\hspace{-0.5mm}1}}**@{},
   <0mm,0mm>*{};<9.0mm,5.5mm>*{^m}**@{},
 <0mm,0mm>*{};<-8mm,-5mm>*{}**@{-},
 <0mm,0mm>*{};<-4.5mm,-5mm>*{}**@{-},
 <0mm,0mm>*{};<-1mm,-5mm>*{\ldots}**@{},
 <0mm,0mm>*{};<4.5mm,-5mm>*{}**@{-},
 <0mm,0mm>*{};<8mm,-5mm>*{}**@{-},
   <0mm,0mm>*{};<-8.5mm,-6.9mm>*{^1}**@{},
   <0mm,0mm>*{};<-5mm,-6.9mm>*{^2}**@{},
   <0mm,0mm>*{};<4.5mm,-6.9mm>*{^{n\hspace{-0.5mm}-\hspace{-0.5mm}1}}**@{},
   <0mm,0mm>*{};<9.0mm,-6.9mm>*{^n}**@{},
 \end{xy}:= \delta
\begin{xy}
 <0mm,0mm>*{\bullet};<0mm,0mm>*{}**@{},
 <0mm,0mm>*{};<-8mm,5mm>*{}**@{-},
 <0mm,0mm>*{};<-4.5mm,5mm>*{}**@{-},
 <0mm,0mm>*{};<-1mm,5mm>*{\ldots}**@{},
 <0mm,0mm>*{};<4.5mm,5mm>*{}**@{-},
 <0mm,0mm>*{};<8mm,5mm>*{}**@{-},
   <0mm,0mm>*{};<-8.5mm,5.5mm>*{^1}**@{},
   <0mm,0mm>*{};<-5mm,5.5mm>*{^2}**@{},
   <0mm,0mm>*{};<4.5mm,5.5mm>*{^{m\hspace{-0.5mm}-\hspace{-0.5mm}1}}**@{},
   <0mm,0mm>*{};<9.0mm,5.5mm>*{^m}**@{},
 <0mm,0mm>*{};<-8mm,-5mm>*{}**@{-},
 <0mm,0mm>*{};<-4.5mm,-5mm>*{}**@{-},
 <0mm,0mm>*{};<-1mm,-5mm>*{\ldots}**@{},
 <0mm,0mm>*{};<4.5mm,-5mm>*{}**@{-},
 <0mm,0mm>*{};<8mm,-5mm>*{}**@{-},
   <0mm,0mm>*{};<-8.5mm,-6.9mm>*{^1}**@{},
   <0mm,0mm>*{};<-5mm,-6.9mm>*{^2}**@{},
   <0mm,0mm>*{};<4.5mm,-6.9mm>*{^{n\hspace{-0.5mm}-\hspace{-0.5mm}1}}**@{},
   <0mm,0mm>*{};<9.0mm,-6.9mm>*{^n}**@{},
 \end{xy}
+
\overset{m-1}{\underset{i=0}{\sum}}
\begin{xy}
 <0mm,0mm>*{\bullet};<0mm,0mm>*{}**@{},
 <0mm,0mm>*{};<-8mm,5mm>*{}**@{-},
 <0mm,0mm>*{};<-3.5mm,5mm>*{}**@{-},
 <0mm,0mm>*{};<-6mm,5mm>*{..}**@{},
 <0mm,0mm>*{};<0mm,5mm>*{}**@{-},
  <0mm,5mm>*{\bullet};
  <0mm,5mm>*{};<0mm,8mm>*{}**@{-},
  <0mm,5mm>*{};<0mm,9mm>*{^{i\hspace{-0.2mm}+\hspace{-0.5mm}1}}**@{},
<0mm,0mm>*{};<8mm,5mm>*{}**@{-},
<0mm,0mm>*{};<3.5mm,5mm>*{}**@{-},
 <0mm,0mm>*{};<6mm,5mm>*{..}**@{},
   <0mm,0mm>*{};<-8.5mm,5.5mm>*{^1}**@{},
   <0mm,0mm>*{};<-4mm,5.5mm>*{^i}**@{},
   <0mm,0mm>*{};<9.0mm,5.5mm>*{^m}**@{},
 <0mm,0mm>*{};<-8mm,-5mm>*{}**@{-},
 <0mm,0mm>*{};<-4.5mm,-5mm>*{}**@{-},
 <0mm,0mm>*{};<-1mm,-5mm>*{\ldots}**@{},
 <0mm,0mm>*{};<4.5mm,-5mm>*{}**@{-},
 <0mm,0mm>*{};<8mm,-5mm>*{}**@{-},
   <0mm,0mm>*{};<-8.5mm,-6.9mm>*{^1}**@{},
   <0mm,0mm>*{};<-5mm,-6.9mm>*{^2}**@{},
   <0mm,0mm>*{};<4.5mm,-6.9mm>*{^{n\hspace{-0.5mm}-\hspace{-0.5mm}1}}**@{},
   <0mm,0mm>*{};<9.0mm,-6.9mm>*{^n}**@{},
 \end{xy}
 - (-1)^{|a|}
\overset{n-1}{\underset{i=0}{\sum}}
 \begin{xy}
 <0mm,0mm>*{\bullet};<0mm,0mm>*{}**@{},
 <0mm,0mm>*{};<-8mm,-5mm>*{}**@{-},
 <0mm,0mm>*{};<-3.5mm,-5mm>*{}**@{-},
 <0mm,0mm>*{};<-6mm,-5mm>*{..}**@{},
 <0mm,0mm>*{};<0mm,-5mm>*{}**@{-},
  <0mm,-5mm>*{\bullet};
  <0mm,-5mm>*{};<0mm,-8mm>*{}**@{-},
  <0mm,-5mm>*{};<0mm,-10mm>*{^{i\hspace{-0.2mm}+\hspace{-0.5mm}1}}**@{},
<0mm,0mm>*{};<8mm,-5mm>*{}**@{-},
<0mm,0mm>*{};<3.5mm,-5mm>*{}**@{-},
 <0mm,0mm>*{};<6mm,-5mm>*{..}**@{},
   <0mm,0mm>*{};<-8.5mm,-6.9mm>*{^1}**@{},
   <0mm,0mm>*{};<-4mm,-6.9mm>*{^i}**@{},
   <0mm,0mm>*{};<9.0mm,-6.9mm>*{^n}**@{},
 <0mm,0mm>*{};<-8mm,5mm>*{}**@{-},
 <0mm,0mm>*{};<-4.5mm,5mm>*{}**@{-},
 <0mm,0mm>*{};<-1mm,5mm>*{\ldots}**@{},
 <0mm,0mm>*{};<4.5mm,5mm>*{}**@{-},
 <0mm,0mm>*{};<8mm,5mm>*{}**@{-},
   <0mm,0mm>*{};<-8.5mm,5.5mm>*{^1}**@{},
   <0mm,0mm>*{};<-5mm,5.5mm>*{^2}**@{},
   <0mm,0mm>*{};<4.5mm,5.5mm>*{^{m\hspace{-0.5mm}-\hspace{-0.5mm}1}}**@{},
   <0mm,0mm>*{};<9.0mm,5.5mm>*{^m}**@{},
 \end{xy}.
 $$
where $\delta$ is the original differential in $\cP$.
The dg properad $(\cP^+, \delta^+)$ is uniquely characterized by the property: there is a 1-1 correspondence between representations
 $$
 \rho: \cP^+ \lon \cE nd_V
 $$
of $(\cP^+, \delta^+)$ in a dg vector space $(V,d)$, and representations of $\cP$ in the same space $V$
but equipped with a deformed differential $d+D$, where $D:=\rho(\begin{xy}
 <0mm,-0.55mm>*{};<0mm,-3mm>*{}**@{-},
 <0mm,0.5mm>*{};<0mm,3mm>*{}**@{-},
 <0mm,0mm>*{\bullet};<0mm,0mm>*{}**@{},
 \end{xy})$.
Clearly any $\cP$-algebra is a $\cP^+$-algebra by letting $D$ act trivially, so that we have a properad map $\cP^+\to \cP$.
Now, slightly abusively, we define $\Der(\cP)$ as the complex of derivations of $\cP^+$ preserving the map $\cP^+\to \cP$. Concretely, in all relevant cases $\cP=\Omega(\cC)$ is the cobar construction of a coaugmented coproperad $\cC$. The definition is then made such that $\Der(\cP)[-1]$ is identified with \eqref{equ:Defdefi} as a complex. On the other hand, if we were using ordinary derivations we would have to modify \eqref{equ:Defdefi} by replacing $\cC$ by the cokernel of the coaugmentation $\overline{\cC}$ on the right-hand side, thus complicating statements of several results.
We assure the reader that this modification is minor and made for technical reasons in the cases we consider, and results about our $\Der(\cP)$ can be easily transcribed into results about the ordinary derivations if necessary.
Note however that $\Der(\cP)$ carries a natural Lie bracket through the commutator.

The deformation complex of a (wheeled) properad $\cP$ is by definition the dg Lie algebra $\Der(\tilde \cP)$ of derivations of a cofibrant resolution $\tilde \cP\stackrel{\sim}{\to}\cP$. It may be identified as a complex with the deformation complex of the identity map $\tilde \cP\to \tilde \cP$ (which controls deformations of $\cP$-algebras) up to a degree shift:
\[
 \Der(\tilde \cP) \cong \Def(\tilde \cP\to \tilde \cP)[1].
\]
Note however that both $\Der(\cP)$ and $\Def(\tilde \cP\to \tilde \cP)$ have natural dg Lie (or $\caL ie_\infty$) algebra structures that are \emph{not} preserved by the above map. Furthermore, there is a quasi-isomorphism of dg Lie algebras
\begin{equation}\label{equ:Defsimpl}
 \Def(\tilde \cP\to \tilde \cP)\to \Def(\tilde \cP\to \cP)
\end{equation}

The zeroth cohomology $H^0(\Der(\tilde \cP))$ is of particular importance. It is a differential graded Lie algebra whose elements act on the space of
$\tilde \cP$-algebra structures on any vector space. We shall see that in the examples we are interested in this dg Lie algebra is very rich, and that it acts non-trivially in general.

\sip
Using the explicit structure of the minimal resolutions of the properads $\LBcd$ and $\LoBcd$ (see \S 2.2 above) we can write down explicit models for the deformation complexes,
\begin{align*}
\Der(\HoLBcd) &=
 \prod_{m,n\geq 1} \left(\HoLBcd(m,n) \otimes \sgn_m^{\ot |c|}\otimes \sgn_n^{\ot |d|}\right)^{\bS_m\times \bS_m}[1+c(1-m)+d(1-n)]
\\
\Der(\HoLoBcd) &= \prod_{m,n\geq 1} \left(\HoLoBcd(m,n)\otimes \sgn_m^{\ot |c|}\otimes \sgn_n^{\ot |d|} \right)^{\bS_m\times \bS_n}[[\hbar]] [1+c(1-m)+d(1-n)]
\end{align*}
Here $\hbar$ is a formal variable of degree $c+d$.
Each of the models on the right has a natural combinatorial interpretation as a graph complex.
For example, $\Der(\HoLBcd)$ may be interpreted as a complex of directed  graphs which have incoming and outgoing legs but have no closed paths of directed edges, for example
$$
\Ga= \resizebox{15mm}{!}{ \xy
(0,0)*{\bu}="d1",
(10,0)*{\bu}="d2",
(-5,-5)*{}="dl",
(5,-5)*{}="dc",
(15,-5)*{}="dr",
(0,10)*{\bu}="u1",
(10,10)*{\bu}="u2",
(5,15)*{}="uc",
(15,15)*{}="ur",
(0,15)*{}="ul",
\ar @{<-} "d1";"d2" <0pt>
\ar @{<-} "d1";"dl" <0pt>
\ar @{<-} "d1";"dc" <0pt>
\ar @{<-} "d2";"dc" <0pt>
\ar @{<-} "d2";"dr" <0pt>
\ar @{<-} "u1";"d1" <0pt>
\ar @{<-} "u1";"d2" <0pt>
\ar @{<-} "u2";"d2" <0pt>
\ar @{<-} "u2";"d1" <0pt>
\ar @{<-} "uc";"u2" <0pt>
\ar @{<-} "ur";"u2" <0pt>
\ar @{<-} "ul";"u1" <0pt>
\endxy} \in \Der(\HoLBcd)
$$
 The value of the differential $\delta$ on an element $\Ga\in \Der(\HoLBcd)$ is obtained by splitting vertices of $\Ga$  and by attaching new corollas at each single external leg of $\Ga$,
$$
\delta \Gamma =
 \delta_{\HoLBcd}\Gamma
 \pm
 \sum\Ba{c}
 \resizebox{9mm}{!}{ \xy
 (0,0)*+{\Ga}="Ga",
(-5,5)*{\bu}="0",
(-8,2)*{}="-1",
(-8,8)*{}="1",
(-5,8)*{}="2",
(-2,8)*{}="3",
\ar @{-} "0";"Ga" <0pt>
\ar @{-} "0";"-1" <0pt>
\ar @{-} "0";"1" <0pt>
\ar @{-} "0";"2" <0pt>
\ar @{-} "0";"3" <0pt>
 \endxy}\Ea
  \pm
 \sum\Ba{c}
\resizebox{9mm}{!}{  \xy
 (0,0)*+{\Ga}="Ga",
(-5,-5)*{\bu}="0",
(-8,-2)*{}="-1",
(-8,-8)*{}="1",
(-5,-2)*{}="2",
(-2,-8)*{}="3",
\ar @{-} "0";"Ga" <0pt>
\ar @{-} "0";"-1" <0pt>
\ar @{-} "0";"1" <0pt>
\ar @{-} "0";"2" <0pt>
\ar @{-} "0";"3" <0pt>
 \endxy}\Ea
 $$
Here $\delta_{\HoLBcd}$ acts on the vertices of $\Ga$ by formula (\ref{LBk_infty}).

\sip

Similarly, $\Der(\HoLoBcd)$ may be interpreted as a complex of $\hbar$-power series of graphs with weighted vertices, for example,
$$
\Ga= \resizebox{15mm}{!}{  \xy
(0,0)*+{_3}*\cir{}="d1",
(10,0)*+{_2}*\cir{}="d2",
(-5,-5)*{}="dl",
(5,-5)*{}="dc",
(15,-5)*{}="dr",
(0,10)*+{_0}*\cir{}="u1",
(10,10)*+{_3}*\cir{}="u2",
(5,15)*{}="uc",
(15,15)*{}="ur",
(0,15)*{}="ul",
\ar @{<-} "d1";"d2" <0pt>
\ar @{<-} "d1";"dl" <0pt>
\ar @{<-} "d1";"dc" <0pt>
\ar @{<-} "d2";"dc" <0pt>
\ar @{<-} "d2";"dr" <0pt>
\ar @{<-} "u1";"d1" <0pt>
\ar @{<-} "u1";"d2" <0pt>
\ar @{<-} "u2";"d2" <0pt>
\ar @{<-} "u2";"d1" <0pt>
\ar @{<-} "uc";"u2" <0pt>
\ar @{<-} "ur";"u2" <0pt>
\ar @{<-} "ul";"u1" <0pt>
\endxy}
\in \Der(\HoLoBcd).
$$
 The value of the differential $\delta$ on an element $\Ga\in \Der(\HoLoBcd)$ consists of three summands,
$$
 \delta \Gamma=
 \delta_{\HoLoBcd}\Gamma
\pm
 \sum \hbar^{p+k-1}\Ba{c}
\resizebox{10mm}{!}{  \xy
  (-3,4)*+{_k},
 (0,0)*+{\Ga}="Ga",
(-7,7)*+{_p}*\cir{}="0",
(-10,4)*{}="-1",
(-10,11)*{}="1",
(-7,11)*{}="2",
(-4,11)*{}="3",
(-5.3,6.8)*-{};(-0.2,1.5)*-{}
**\crv{(0.4,5)};
(-6.9,5.1)*-{};(-1.4,0.3)*-{}
**\crv{(-5,-1.5)};
(-6.0,5.4)*-{};(-1.4,0.7)*-{}
**\crv{(-4,-0.5)};
%
\ar @{-} "0";"-1" <0pt>
\ar @{-} "0";"1" <0pt>
\ar @{-} "0";"2" <0pt>
\ar @{-} "0";"3" <0pt>
 \endxy}\Ea
  \pm
 \sum \hbar^{p+k-1}\Ba{c}
 \resizebox{10mm}{!}{  \xy
  (-3,-4)*+{_k},
 (0,0)*+{\Ga}="Ga",
(-7,-7)*+{_p}*\cir{}="0",
(-10,-4)*{}="-1",
(-10,-11)*{}="1",
(-7,-11)*{}="2",
(-4,-11)*{}="3",
(-5.3,-6.8)*-{};(-0.2,-1.5)*-{}
**\crv{(0.4,-5)};
(-6.9,-5.1)*-{};(-1.4,-0.3)*-{}
**\crv{(-5,1.5)};
(-6.0,-5.4)*-{};(-1.4,-0.7)*-{}
**\crv{(-4,-0.5)};
%
\ar @{-} "0";"-1" <0pt>
\ar @{-} "0";"1" <0pt>
\ar @{-} "0";"2" <0pt>
\ar @{-} "0";"3" <0pt>
 \endxy}\Ea
$$
 where the first summand comes from the action of the $\HoLoBcd$-differential
 on internal vertices of $\Ga$ by formula (\ref{2: d on Lie inv infty}), and in
 the two terms on the right one sums over all ways of attaching a new vertex to some subset of the incoming or outgoing legs ($k$ many), and sums over all possible decorations $p$ of the added vertex, with an appropriate power of $\hbar$ as prefactor. Note that the power of $\hbar$ counts the number of loops added to the graph, if we count a vertex decorated by $p$ as contributing $p$ loops.

 \sip

The Lie bracket is combinatorially obtained by inserting graphs into vertices of another.

\sip

The cohomology of all these graph complexes is hard to compute. We may however simplify the computation of three of the above four complexes by using formula (\ref{equ:Defsimpl}) and equivalently study instead the following much ``smaller" complexes,

\begin{align*}
\Def(\HoLBcd\rar \LBcd) &=
 \prod_{m,n\geq 1} \left(\LBcd(m,n) \otimes \sgn_m^{\ot |c|}\otimes \sgn_n^{\ot |d|}\right)^{\bS_m\times \bS_m}[c(1-m)+d(1-n)]
\\
\Def(\HoLoBcd\rar \LoBcd) &= \prod_{m,n\geq 1} \left(\LoBcd(m,n)\otimes \sgn_m^{\ot |c|}\otimes \sgn_n^{\ot |d|} \right)^{\bS_m\times \bS_n} [[\hbar]][c(1-m)+d(1-n)]
\end{align*}

Note however that in passing from $\Der(\dots)$ to the (quasi-isomorphic) simpler complexes $\Def(\dots)$ above we lose the dg Lie algebra structure, or rather  there is a different Lie algebra structure on the above complexes.
The above complexes may again be interpreted as graph complexes. For example $\Def(\HoLBcd\to \LBcd)$ consists of oriented trivalent graphs with incoming and outgoing legs, modulo the Jacobi and Drinfeld five term relations. The differential is obtained by attaching a trivalent vertex at one external leg in all possible ways.

\subsection{Complete variants}
It is more convenient for our purposes to consider (genus-)completed versions of the deformation complexes of the previous subsection. In particular, the the genus filtration endows the properads $\wHoLBcd$ and $\wHoLoBcd$  with complete topologies, and we define the complexes of \emph{continuous} derivations $\Der(\wHoLBcd)$ and $\Der(\wHoLoBcd)$.

\sip

More concretely, the difference of the complete and incomplete versions is as follows:
\begin{itemize}
 \item Elements of $\Der(\HoLBcd)$  of a fixed cohomological degree can be identified with possible infinite series of graphs. However, in each fixed arity, the series must be a finite linear combination.
 On the other hand, elements of $\Der(\wHoLBcd)$  in a fixed degree are simply all series of graphs.


 \item Similarly, elements of $\Der(\HoLoBcd)$ may be understood as power series in $\hbar$ with coefficients in the
 series of graphs with the same finiteness condition as before.
 On the other hand elements of $\Der(\wHoLoBcd)$ are power series in $\hbar$ with coefficients arbitrary series of graphs.
\end{itemize}

Below we will always work with the complete versions of our properads and deformation/derivation complexes.

\subsection{\bf A map from the  graph complex $\GCor_{c+d+1}$ to $\Der(\wHoLBcd)$}\label{3: subsect on GCor_3 to Der(LieB)} A derivation $D$ of any free properad $\cP$
is uniquely determined by its values, $D(e)\in \cP$, on the generators $e$ of $\cP$.

\sip

There is a natural right action of the dg Lie algebra $\mathsf{fGC}^{or}_{c+d+1}$ on the genus completed dg properad $\wHoLBcd$ by properadic derivations (cf.\ \cite{CMW}), i.e.\ there is a canonical morphism
of dg Lie algebras,
\Beq\label{2: Morhism F from GC_3^or}
\Ba{rccc}
 F\colon & \mathsf{fGC}^{or}_{c+d+1} &\to & \Der(\wHoLBcd)\\
         &   \Ga & \to & F(\Ga)
         \Ea
\Eeq
with  values of the derivation $F(\Ga)$
on the generators of the (genus) completed properad  $\hLieBi_\infty$
given explicitly by
\Beq \label{equ:def GC action 1}
\left(\Ba{c}\resizebox{12mm}{!}{\begin{xy}
 <0mm,0mm>*{\circ};<0mm,0mm>*{}**@{},
 <-0.6mm,0.44mm>*{};<-8mm,5mm>*{}**@{-},
 <-0.4mm,0.7mm>*{};<-4.5mm,5mm>*{}**@{-},
 <0mm,0mm>*{};<-1mm,5mm>*{\ldots}**@{},
 <0.4mm,0.7mm>*{};<4.5mm,5mm>*{}**@{-},
 <0.6mm,0.44mm>*{};<8mm,5mm>*{}**@{-},
   <0mm,0mm>*{};<-8.5mm,5.5mm>*{^1}**@{},
   <0mm,0mm>*{};<-5mm,5.5mm>*{^2}**@{},
   <0mm,0mm>*{};<9.0mm,5.5mm>*{^m}**@{},
 <-0.6mm,-0.44mm>*{};<-8mm,-5mm>*{}**@{-},
 <-0.4mm,-0.7mm>*{};<-4.5mm,-5mm>*{}**@{-},
 <0mm,0mm>*{};<-1mm,-5mm>*{\ldots}**@{},
 <0.4mm,-0.7mm>*{};<4.5mm,-5mm>*{}**@{-},
 <0.6mm,-0.44mm>*{};<8mm,-5mm>*{}**@{-},
   <0mm,0mm>*{};<-8.5mm,-6.9mm>*{^1}**@{},
   <0mm,0mm>*{};<-5mm,-6.9mm>*{^2}**@{},
   <0mm,0mm>*{};<9.0mm,-6.9mm>*{^n}**@{},
 \end{xy}}\Ea\right)\cdot F(\Ga)
=
 \sum_{s:[n]\rar V(\Ga)\atop \hat{s}:[m]\rar V(\Ga)}  \Ba{c}\resizebox{9mm}{!}  {\xy
 (-6,7)*{^1},
(-3,7)*{^2},
(2.5,7)*{},
(7,7)*{^m},
(-3,-8)*{_2},
(3,-6)*{},
(7,-8)*{_n},
(-6,-8)*{_1},
(0,4.5)*+{...},
(0,-4.5)*+{...},
(0,0)*+{\Ga}="o",
(-6,6)*{}="1",
(-3,6)*{}="2",
(3,6)*{}="3",
(6,6)*{}="4",
(-3,-6)*{}="5",
(3,-6)*{}="6",
(6,-6)*{}="7",
(-6,-6)*{}="8",
\ar @{-} "o";"1" <0pt>
\ar @{-} "o";"2" <0pt>
\ar @{-} "o";"3" <0pt>
\ar @{-} "o";"4" <0pt>
\ar @{-} "o";"5" <0pt>
\ar @{-} "o";"6" <0pt>
\ar @{-} "o";"7" <0pt>
\ar @{-} "o";"8" <0pt>
\endxy}\Ea
\Eeq
 where the sum is taken over all ways of attaching the incoming and outgoing legs to the graph $\Ga$, and we set to zero every resulting graph if it contains a vertex with valency $<3$ or
   with no at least one incoming  or at least one outgoing edge.

   \sip

   Let us first check that the map $f$ has degree zero. If $\Ga\in \mathsf{fGC}_{c+d+1}^{or}$ has $p$ vertices and $l$ edges, then $|\Ga|=(c+d+1)(p-1) - (c+d)l$ so that the total degree of the l.h.s.\ in
   (\ref{equ:def GC action 1}) equals
   $$
   (c+d+1)(p-1) - (c+d)l + 1 +c(1-m)+d(1-n)= (c+d+1)p - (c+d)l-cm-dn.
   $$
   On the other hand, each summand on the r.h.s.\ gives us a graph in $ \LieBi_\infty$
   with $p$ vertices, and each vertex $v$ has $|v|_{out} + \# \hat{s}^{-1}(v)$ output edges and $|v|_{in}+ \# {s}^{-1}(v)$ input edges,
 where  $|v|_{in}$ (resp.\ $|v|_{out}$) counts the number of outputs (resp.\ inputs) of $v$ in $\Ga$. The degree of such a graph is equal to
   \Beqrn
   \sum_{v\in V(\Ga)}
    \left(1+c(1-|v|_{out} - \# \hat{s}^{-1}(v)) + d(1-|v|_{in}- \# s^{-1}(v))\right)&=&p +c(p-l-m) +d(p-l-n)\\
    &=& (c+d+1)p- (c+d)l -cm-dn
   \Eeqrn
as required. There is an implicit rule of signs in-built into formula (\ref{equ:def GC action 1}) which is completely analogous to the one defined in \S 7 of \cite{MaVo}.

 \sip

Consider for a moment $\wHoLBcd$ as a {\em non-differential} \, (completed) free properad, and let $\Der_{non-d}(\wHoLBcd)$ stand for its Lie algebra of derivations
which can be identified with the vector space
$$
\prod_{m,n} \Hom_{\bS_m^{op}\times \bS_n}\left(\K\langle {\mathfrak C_n^m} \rangle,\wHoLBcd(m,n)\right)=\prod_{m,n\geq 1} \left(\LBcd(m,n) \otimes \sgn_m^{\ot |c|}\otimes \sgn_n^{\ot |d|}\right)^{\bS_m\times \bS_m}[c(1-m)+d(1-n)]
$$
as any derivation is uniquely determined by its values on the generators
$$
{\mathfrak C_n^m}:=\Ba{c}\resizebox{12mm}{!}{\begin{xy}
 <0mm,0mm>*{\circ};<0mm,0mm>*{}**@{},
 <-0.6mm,0.44mm>*{};<-8mm,5mm>*{}**@{-},
 <-0.4mm,0.7mm>*{};<-4.5mm,5mm>*{}**@{-},
 <0mm,0mm>*{};<-1mm,5mm>*{\ldots}**@{},
 <0.4mm,0.7mm>*{};<4.5mm,5mm>*{}**@{-},
 <0.6mm,0.44mm>*{};<8mm,5mm>*{}**@{-},
   <0mm,0mm>*{};<-8.5mm,5.5mm>*{^1}**@{},
   <0mm,0mm>*{};<-5mm,5.5mm>*{^2}**@{},
   <0mm,0mm>*{};<9.0mm,5.5mm>*{^m}**@{},
 <-0.6mm,-0.44mm>*{};<-8mm,-5mm>*{}**@{-},
 <-0.4mm,-0.7mm>*{};<-4.5mm,-5mm>*{}**@{-},
 <0mm,0mm>*{};<-1mm,-5mm>*{\ldots}**@{},
 <0.4mm,-0.7mm>*{};<4.5mm,-5mm>*{}**@{-},
 <0.6mm,-0.44mm>*{};<8mm,-5mm>*{}**@{-},
   <0mm,0mm>*{};<-8.5mm,-6.9mm>*{^1}**@{},
   <0mm,0mm>*{};<-5mm,-6.9mm>*{^2}**@{},
   <0mm,0mm>*{};<9.0mm,-6.9mm>*{^n}**@{},
 \end{xy}}\Ea
$$
which can be chosen arbitrary.
It is immediate from the definition of the Lie algebra structure $[\ ,\ ]$ in $\GCor_{c+d+1}$
that the map
$$
F: \mathsf{fGC}_{c+d+1}^{or} \lon \Der_{non-d}(\hLieBi_\infty)
$$
given by the formula (\ref{equ:def GC action 1})
respects the Lie brackets,
$$
\left(\mathfrak C_n^m \cdot F(\Ga_1)\right)\cdot F(\Ga_2) - (-1)^{|\Ga_1||\Ga_2|}
\left(\mathfrak C_n^m \cdot F(\Ga_2)\right)\cdot F(\Ga_1) = \mathfrak C_n^m \cdot F([\Ga_1,\Ga_2])
$$
for any $\Ga_1,\Ga_2\in \GCor_{c+d+1}$. This result implies that any Maurer-Cartan element $\Ga$ in the Lie algebra
$(\mathsf{fGC}^{or}_{c+d+1}, [\ ,\ ])$ gives rise to a continuous differential
$$
d_\Ga: \mathfrak C_n^m \lon \mathfrak C_n^m \cdot F(\Ga)
$$
in the properad $\wHoLBcd$. A remarkable (and almost obvious) fact is that the Maurer-Cartan element
$$
\Ga= \xy
  (0,-3)*{\bullet}="a",
 (0,3)*{\bu}="b",
 \ar @{->} "a";"b" <0pt> \endxy
$$
induces the standard differential (\ref{LBk_infty}) in $\wHoLBcd$. This implies that the morphism (\ref{2: Morhism F from GC_3^or}) induces (by changing the right action into a left action via  a standard sign factor)  a map of {\em dg}\, Lie algebras
\Beq\label{2: F final from GC_3^or to DerLieBi}
 F\colon \sG\sC_3^{c+d+1}\to \Der(\wHoLBcd)\,
\Eeq
 which is proven below to be a quasi-isomorphism (up to one class).

\subsubsection{\bf Remark}
Interpreting the right hand side in (\ref{2: F final from GC_3^or to DerLieBi})   as a graph complex itself (see section  {\S \ref{2 sec:defcomplexes}}), we see that the map $F$ sends a graph $\Gamma\in \GC_{c+d+1}^{or}$ to the series of graphs
\[
\sum_{m,n\geq 1}
 \sum_{s:[n]\rar V(\Ga)\atop \hat{s}:[m]\rar V(\Ga)}
    \overbrace{
 \underbrace{\Ba{c}\resizebox{10mm}{!}  { \xy
(0,4.5)*+{...},
(0,-4.5)*+{...},
(0,0)*+{\Ga}="o",
(-5,6)*{}="1",
(-3,6)*{}="2",
(3,6)*{}="3",
(5,6)*{}="4",
(-3,-6)*{}="5",
(3,-6)*{}="6",
(5,-6)*{}="7",
(-5,-6)*{}="8",
\ar @{-} "o";"1" <0pt>
\ar @{-} "o";"2" <0pt>
\ar @{-} "o";"3" <0pt>
\ar @{-} "o";"4" <0pt>
\ar @{-} "o";"5" <0pt>
\ar @{-} "o";"6" <0pt>
\ar @{-} "o";"7" <0pt>
\ar @{-} "o";"8" <0pt>
\endxy}\Ea
 }_{n\times}
 }^{m\times}
%
%
%
 %
\]
where the  second summation symbol has exactly the same meaning as in (\ref{equ:def GC action 1}).

\subsection{\bf A map from the  graph complex $\GCor_{c+d+1}[[\hbar]]$ to {$\Der(\widehat{\HoLoBcd})$}}
There is a natural right action of the {\em non-differential}\, Lie algebra $(\GC_{c+d+1}^{or}[[\hbar]], [\ ,\ ])$  on  the {\em non-differential}\,
free operad $\HoLoBcd$ by continuous derivations, that is, there is a continuous morphism of
topological Lie algebras,
$$
F^\diamond: \left( \mathsf{fGC}^{or}_{c+d+1}[[\hbar]], [\ ,\ ]\right) \lon \Der(\wHoLoBcd)
$$
For any monomial $\hbar^k\Ga\in \GC_{c+d+1}^{or}[[\hbar]]$ the value of the associated derivation $F^\diamond(\hbar^k \Ga)$
on  the generators of   $\wHoLoBcd$
is given, by definition,  by
\Beq \label{equ:def GC action 2}
\left(
\Ba{c}\resizebox{14mm}{!}{\xy
(-9,-6)*{};
(0,0)*+{a}*\cir{}
**\dir{-};
(-5,-6)*{};
(0,0)*+{a}*\cir{}
**\dir{-};
(9,-6)*{};
(0,0)*+{a}*\cir{}
**\dir{-};
(5,-6)*{};
(0,0)*+{a}*\cir{}
**\dir{-};
(0,-6)*{\ldots};
(-10,-8)*{_1};
(-6,-8)*{_2};
(10,-8)*{_n};
(-9,6)*{};
(0,0)*+{a}*\cir{}
**\dir{-};
(-5,6)*{};
(0,0)*+{a}*\cir{}
**\dir{-};
(9,6)*{};
(0,0)*+{a}*\cir{}
**\dir{-};
(5,6)*{};
(0,0)*+{a}*\cir{}
**\dir{-};
(0,6)*{\ldots};
(-10,8)*{_1};
(-6,8)*{_2};
(10,8)*{_m};
\endxy}\Ea
\right)\cdot F^(\hbar^k\Ga)
:=
\left\{\Ba{ll} \displaystyle
 \sum_{s:[n]\rar V(\Ga)\atop \hat{s}:[m]\rar V(\Ga)}
 \sum_{a=k+\sum_{v\in V(\Ga)} a_v\atop a_v\geq 0}
  \Ba{c}\resizebox{9mm}{!}  {\xy
 (-6,7)*{^1},
(-3,7)*{^2},
(2.5,7)*{},
(7,7)*{^m},
(-3,-8)*{_2},
(3,-6)*{},
(7,-8)*{_n},
(-6,-8)*{_1},
(0,4.5)*+{...},
(0,-4.5)*+{...},
(0,0)*+{\Ga}="o",
(-6,6)*{}="1",
(-3,6)*{}="2",
(3,6)*{}="3",
(6,6)*{}="4",
(-3,-6)*{}="5",
(3,-6)*{}="6",
(6,-6)*{}="7",
(-6,-6)*{}="8",
\ar @{-} "o";"1" <0pt>
\ar @{-} "o";"2" <0pt>
\ar @{-} "o";"3" <0pt>
\ar @{-} "o";"4" <0pt>
\ar @{-} "o";"5" <0pt>
\ar @{-} "o";"6" <0pt>
\ar @{-} "o";"7" <0pt>
\ar @{-} "o";"8" <0pt>
\endxy}\Ea
  & \ \ \ \ \ \mbox{if}\ \ k\leq a\\
 0 & \ \ \ \ \  \mbox{if}\ \  k>a,
 \Ea
 \right.
\Eeq
where the first sum is taking over all ways to attach $m$ output legs and $n$ input legs to the vertices
of the graph $\Ga$, and the second sum is taken over all  ways to decorate the vertices of $\Ga$ with non-negative integers $a_1,\ldots,a_{\# V(\Ga)}$ such they sum to $a-k$; moreover, we set a graph on the r.h.s.\ to zero is  if there is at least one vertex $v$ with the number $n_v$ of incoming edges equal to zero, or the number $m_v$ of outgoing number equal to zero, or if the condition $n_v+m_v+a_v\geq 3$ is violated. There is an implicit rule of signs in formula
(\ref{equ:def GC action 2}) which is identical to the one in the subsection above.

\sip

It is easy to check that $f^\diamond$ has degree zero and respects Lie brackets. Therefore, it sends any
Maurer-Cartan element in  $( \mathsf{fGC}^{or}_3[[\hbar]], [\ ,\ ])$ into a differential in the free prop
$\wHoLoBcd$. It is again almost immediate to see that the differential induced by the
 Maurer-Cartan element (\ref{2: Phi_hbar MC element}) is precisely the one given in (\ref{2: d on Lie inv infty}), i.e.\ the one which makes $\wHoLoBcd$ into a minimal resolution of $\wLoBcd$.
Therefore we conclude that there is a morphism of {\em dg}\, Lie algebras
\Beq\label{3: F from GCor[[hbar]] to Der LoBinfty}
f^\diamond: \left(\GCor_3[[\hbar]], \delta_\hbar=[\Phi_\hbar,\ ]\right) \lon \left(\Der(\hLoB_\infty), d:=[\delta,\ ]\right).
\Eeq
We shall prove below that this map is almost   a quasi-isomorphism.

\subsubsection{\bf Remark}
Interpreting $\Der(\wHoLoBcd)$   as a graph complex  (see section  {\S \ref{2 sec:defcomplexes}}), we can reformulate  the map $F^\diamond$ as the one which sends a monomial $\hbar^k\Gamma\in \GC_{c+d+1}^{or}[[\hbar]]$ to the series of graphs
\[
\sum_{m,n\geq 1}
 \sum_{s:[n]\rar V(\Ga)\atop \hat{s}:[m]\rar V(\Ga)}
 \sum_{a=k+\sum_{v\in V(\Ga)} a_v\atop a_v\geq 0}
    \overbrace{
 \underbrace{\Ba{c}\resizebox{10mm}{!}  { \xy
(0,4.5)*+{...},
(0,-4.5)*+{...},
(0,0)*+{\Ga}="o",
(-5,6)*{}="1",
(-3,6)*{}="2",
(3,6)*{}="3",
(5,6)*{}="4",
(-3,-6)*{}="5",
(3,-6)*{}="6",
(5,-6)*{}="7",
(-5,-6)*{}="8",
\ar @{-} "o";"1" <0pt>
\ar @{-} "o";"2" <0pt>
\ar @{-} "o";"3" <0pt>
\ar @{-} "o";"4" <0pt>
\ar @{-} "o";"5" <0pt>
\ar @{-} "o";"6" <0pt>
\ar @{-} "o";"7" <0pt>
\ar @{-} "o";"8" <0pt>
\endxy}\Ea
 }_{n\times}
 }^{m\times}
\]
where the  second and third summation symbols have exactly the same meaning as in (\ref{equ:def GC action 2}). We shall use this fact below.

\bip

{
\Large
\section{\bf Computations of the cohomology of deformation complexes}
}

\mip

In this section we compute the cohomology of several of the deformation complexes, show Theorems {\ref{thm:Fqiso}} and {\ref{thm:Fhbarqiso}}, and discuss their concrete applications.

\mip

\subsection{The proof of Theorem {\ref{thm:Fqiso}}} \label{app:defproof1}
Let us recall the definition of the graph complex $\hGCor_d$ from \cite[section 3.3]{Wi2}.
The elements of $\hGCor_d$ are $\K$-linear series in directed acyclic graphs with outgoing legs such that all vertices are at least bivalent, and such that there are no bivalent vertices with one incoming and one outgoing edge.\footnote{The last condition is again not present on \cite{Wi2}, but it does not change the cohomology.} We set to zero graphs containing vertices without outgoing edges. Here is an example graph:
$$
 \Ba{c}\resizebox{6mm}{!} 
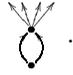
{\xy
 (0,0)*{\bu}="o",
(-5,6)*{}="d1",
(-2,6)*{}="d2",
(2,6)*{}="d3",
(5,6)*{}="d4",
\ar @{->} "o";"d1" <0pt>
\ar @{->} "o";"d2" <0pt>
\ar @{->} "o";"d3" <0pt>
\ar @{->} "o";"d4" <0pt>
   \ar@/^0.6pc/(0,-8)*{\bullet};(0,0)*{\bullet}
   \ar@/^{-0.6pc}/(0,-8)*{\bullet};(0,0)*{\bullet}
 \endxy}
 \Ea
%
.$$
The degrees are computed just as for graphs occurring in $\GCor_d$, with the external legs considered to be of degree 0. For the description of the differential we refer the reader to \cite[section 3.3]{Wi2}.

\sip

There is a map $\Psi: \GCor_d\to \hGCor_d$ sending a graph $\Gamma$ to the linear combination
\begin{equation}\label{equ:hairymap}
\Gamma \mapsto
\sum_{j=1}^\infty
\overbrace{
  \Ba{c}\resizebox{9mm}{!}  {\xy
(0,5.2)*+{...},
(0,0)*+{\Ga}="o",
(-3,7)*{}="5",
(3,7)*{}="6",
(5,7)*{}="7",
(-5,7)*{}="8",
\ar @{->} "o";"5" <0pt>
\ar @{->} "o";"6" <0pt>
\ar @{->} "o";"7" <0pt>
\ar @{->} "o";"8" <0pt>
\endxy}\Ea
 }^{j\times}
\end{equation}
where the picture on the right means that one should sum over all ways of connecting $j$ outgoing edges to the graph $\Gamma$. Graphs for which there remain vertices with no outgoing edge are identified with $0$.

The following proposition has been shown in loc. cit.

\subsubsection{\bf Proposition}[Proposition 3 of \cite{Wi2}]\label{prop:GChGC}
{\em The map $\Psi: \GCor_d \to \hGCor_d$ is a quasi-isomorphism up to the class in $H(\hGCor_d)$ represented by the graph cocycle
\begin{equation}
\label{equ:singleclass}
 \sum_{j\geq 2}
 (j-1)
 \overbrace{
  \Ba{c}\resizebox{9mm}{!}  {\xy
(0,5.2)*+{...},
(0,0)*{\bu}="o",
(-3,7)*{}="5",
(3,7)*{}="6",
(5,7)*{}="7",
(-5,7)*{}="8",
\ar @{->} "o";"5" <0pt>
\ar @{->} "o";"6" <0pt>
\ar @{->} "o";"7" <0pt>
\ar @{->} "o";"8" <0pt>
\endxy}\Ea
 }^{j\times}.
\end{equation}
}

\mip

There is a map $G:\hGCor_{c+d+1} \to \Der(\wHoLBcd)$ sending a graph $\overbrace{
  \Ba{c}\resizebox{9mm}{!}  {\xy
(0,5.2)*+{...},
(0,0)*+{\Ga}="o",
(-3,7)*{}="5",
(3,7)*{}="6",
(5,7)*{}="7",
(-5,7)*{}="8",
\ar @{->} "o";"5" <0pt>
\ar @{->} "o";"6" <0pt>
\ar @{->} "o";"7" <0pt>
\ar @{->} "o";"8" <0pt>
\endxy}\Ea}^{m\times} \in \hGCor_{c+d+1}$ to the series
\Beq 
G( \Gamma )=
 \sum_{n}
    \overbrace{
 \underbrace{ \Ba{c}\resizebox{9mm}{!}  {\xy
(0,4.9)*+{...},
(0,-4.9)*+{...},
(0,0)*+{\Ga}="o",
(-5,6)*{}="1",
(-3,6)*{}="2",
(3,6)*{}="3",
(5,6)*{}="4",
(-3,-6)*{}="5",
(3,-6)*{}="6",
(5,-6)*{}="7",
(-5,-6)*{}="8",
\ar @{->} "o";"1" <0pt>
\ar @{->} "o";"2" <0pt>
\ar @{->} "o";"3" <0pt>
\ar @{->} "o";"4" <0pt>
\ar @{<-} "o";"5" <0pt>
\ar @{<-} "o";"6" <0pt>
\ar @{<-} "o";"7" <0pt>
\ar @{<-} "o";"8" <0pt>
\endxy}\Ea
 }_{n\times}
 }^{m\times}
\Eeq

The map $F : \GC_{c+d+1}^{or}\to \Der(\wHoLBcd)$ from Theorem {\ref{thm:Fqiso}} factors through the map $G$ above, i.~e., it can be written as the composition
\[
 \GC_{c+d+1}^{or}\stackrel{\Psi}{\to} \hGCor_{c+d+1} \stackrel{G}{\to} \Der(\wHoLBcd).
\]
In view of Proposition {\ref{prop:GChGC}}, Theorem {\ref{thm:Fqiso}} hence follows immediately from the following result.

\subsubsection{\bf Proposition}\label{prop:Gqiso}
{\em  The map $G:\hGCor_{c+d+1} \to \Der(\wHoLBcd)$ is a quasi-isomorphism.
}

\begin{proof}
 For a graph in $\Der(\wHoLBcd)$ we will call its \emph{skeleton} the graph obtained in the following way:
 \begin{enumerate}
  \item Remove all input legs and recursively remove all valence 1 vertices created.
  \item Remove valence 2 vertices with one incoming and one outgoing edge and connect the two edges.
 \end{enumerate}
An example of a graph and its skeleton the following
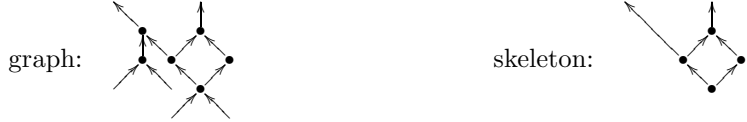
\begin{align*}
 \text{graph: }&
 \Ba{c}
\resizebox{17mm}{!}
{ \xy
(0,-10)*{\bu}="u",
(-5,-5)*{\bu}="L",
 (5,-5)*{\bu}="R",
(0,0)*{\bu}="d",
(-5,-15)*{}="u1",
(5,-15)*{}="u2",
(0,5)*{}="d1",
(-10,0)*{\bu}="b",
(-10,-5)*{\bu}="a",
(-15,-10)*{}="a1",
(-5,-10)*{}="a2",
(-15,5)*{}="b1",
\ar @{->} "u1";"u" <0pt>
\ar @{->} "u2";"u" <0pt>
\ar @{->} "u";"L" <0pt>
\ar @{->} "u";"R" <0pt>
\ar @{->} "L";"d" <0pt>
\ar @{->} "R";"d" <0pt>
\ar @{->} "d";"d1" <0pt>
\ar @{->} "L";"b" <0pt>
\ar @{->} "a1";"a" <0pt>
\ar @{->} "a2";"a" <0pt>
\ar @{->} "b";"b1" <0pt>
\ar @{->} "a";"b" <0pt>
\endxy}
\Ea
 &
  \text{skeleton: }&
   \Ba{c}
\resizebox{17mm}{!}
{ \xy
(0,-10)*{\bu}="u",
(-5,-5)*{\bu}="L",
 (5,-5)*{\bu}="R",
(0,0)*{\bu}="d",
(-5,-15)*{}="u1",
(5,-15)*{}="u2",
(0,5)*{}="d1",
(-15,-10)*{}="a1",
(-5,-10)*{}="a2",
(-15,5)*{}="b1",
\ar @{->} "u";"L" <0pt>
\ar @{->} "u";"R" <0pt>
\ar @{->} "L";"d" <0pt>
\ar @{->} "R";"d" <0pt>
\ar @{->} "d";"d1" <0pt>
\ar @{->} "L";"b1" <0pt>
\endxy}
\Ea
\end{align*}
We put a filtration on $\Der(\wHoLBcd)$ by the total number of vertices in the skeleton. Let $\gr\Der(\wHoLBcd)$ be the associated graded.
Note that for elements in the image of some graph $\Gamma\in \hGCor_{c+d+1}$ under $G$ the skeleton is just the graph $\Gamma$, and hence there is a map of complexes $\hGCor_{c+d+1} \to \gr\Der(\wHoLBcd)$, where we consider the left hand side with zero differential.
We claim that the induced map $\hGCor_{c+d+1} \to H(\gr\Der(\wHoLBcd))$ is an isomorphism. From this claim the Proposition follows immediately by a standard spectral sequence argument.

\sip

The differential on $\gr\Der(\wHoLBcd)$ does not change the skeleton. Hence the complex $\gr\Der(\wHoLBcd)$ splits into a direct product of complexes, say $\tilde C_\gamma$, one for each skeleton $\gamma$
\[
 \gr\Der(\wHoLBcd) = \prod_\gamma \tilde C_\gamma .
\]
Furthermore, each skeleton represents an automorphism class of graphs, and we may write
\[
 \tilde C_\gamma = C_{\tilde \gamma}^{\Aut_\gamma}
\]
where $C_{\tilde \gamma}$ is an appropriately defined complex for one representative $\tilde \gamma$ of the isomorphism class $\gamma$ and $\Aut_\gamma$ is the automorphism group associated to the skeleton. In other words, the $\tilde \gamma$ now has distinguishable vertices and edges.
More concretely, the complex $C_{\tilde \gamma}$ is the complex of $\K$-linear series of graphs obtained from $\tilde \gamma$ by
\begin{enumerate}
 \item Adding some bivalent vertices with one input and one output on edges. We call these vertices ``edge vertices''.
 \item Attaching input forests at the vertices, such that all vertices are at least trivalent and have at least one input and one output.
 We call the forest attached to a vertex the forest of that vertex.
\end{enumerate}
An example is the following:
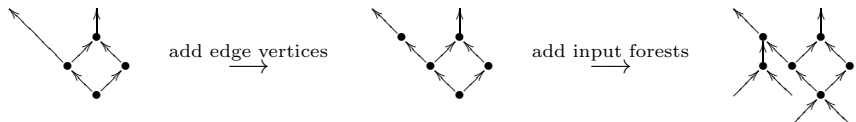
\[
  \Ba{c}
\resizebox{17mm}{!}
{ \xy
(0,-10)*{\bu}="u",
(-5,-5)*{\bu}="L",
 (5,-5)*{\bu}="R",
(0,0)*{\bu}="d",
(-5,-15)*{}="u1",
(5,-15)*{}="u2",
(0,5)*{}="d1",
(-15,-10)*{}="a1",
(-5,-10)*{}="a2",
(-15,5)*{}="b1",
\ar @{->} "u";"L" <0pt>
\ar @{->} "u";"R" <0pt>
\ar @{->} "L";"d" <0pt>
\ar @{->} "R";"d" <0pt>
\ar @{->} "d";"d1" <0pt>
\ar @{->} "L";"b1" <0pt>
\endxy}
\Ea
\ \ \stackrel{\text{add edge vertices}}{\longrightarrow} \ \
  \Ba{c}
\resizebox{17mm}{!}
{ \xy
(0,-10)*{\bu}="u",
(-5,-5)*{\bu}="L",
 (5,-5)*{\bu}="R",
(0,0)*{\bu}="d",
(-5,-15)*{}="u1",
(5,-15)*{}="u2",
(0,5)*{}="d1",
(-10,0)*{\bu}="b",
(-15,-10)*{}="a1",
(-5,-10)*{}="a2",
(-15,5)*{}="b1",
\ar @{->} "u";"L" <0pt>
\ar @{->} "u";"R" <0pt>
\ar @{->} "L";"d" <0pt>
\ar @{->} "R";"d" <0pt>
\ar @{->} "d";"d1" <0pt>
\ar @{->} "L";"b" <0pt>
\ar @{->} "b";"b1" <0pt>
\endxy}
\Ea
 \ \ \stackrel{\text{add input forests}}{\longrightarrow}\ \
 \Ba{c}
\resizebox{17mm}{!}
{ \xy
(0,-10)*{\bu}="u",
(-5,-5)*{\bu}="L",
 (5,-5)*{\bu}="R",
(0,0)*{\bu}="d",
(-5,-15)*{}="u1",
(5,-15)*{}="u2",
(0,5)*{}="d1",
(-10,0)*{\bu}="b",
(-10,-5)*{\bu}="a",
(-15,-10)*{}="a1",
(-5,-10)*{}="a2",
(-15,5)*{}="b1",
\ar @{->} "u1";"u" <0pt>
\ar @{->} "u2";"u" <0pt>
\ar @{->} "u";"L" <0pt>
\ar @{->} "u";"R" <0pt>
\ar @{->} "L";"d" <0pt>
\ar @{->} "R";"d" <0pt>
\ar @{->} "d";"d1" <0pt>
\ar @{->} "L";"b" <0pt>
\ar @{->} "a1";"a" <0pt>
\ar @{->} "a2";"a" <0pt>
\ar @{->} "b";"b1" <0pt>
\ar @{->} "a";"b" <0pt>
\endxy}
\Ea
\]
We next put another filtration on $C_{\tilde \gamma}$ by the number of edge vertices added in the first step above and consider the associated graded $\gr C_{\tilde \gamma}$.
Note that the differential of $\gr C_{\tilde \gamma}$ acts on each of the forests attached to the vertices separately and hence the complex splits into a (completed) tensor product of complexes, one for each such vertex. Let us call the complex made from the possible forests at the vertex $v$ the forest complex at that vertex.
By the same argument showing that the cohomology of a free $\caL ie_\infty$ algebra generated by a single generator is two dimensional, we find that the
forest complex at $v$ has either one or two dimensional cohomology. If the vertex $v$ has no incoming edge then the cohomology of the forest complex is one dimensional, the class being represented by the forests
\[
\sum_{j\geq 1}
\underbrace{
 \Ba{c}\resizebox{9mm}{!}  {\xy
(0,-4.5)*+{...},
(0,2)*{_v},
(0,0)*{\bu}="o",
(-5,-5)*{}="1",
(-3,-5)*{}="2",
(3,-5)*{}="3",
(5,-5)*{}="4",
\ar @{<-} "o";"1" <0pt>
\ar @{<-} "o";"2" <0pt>
\ar @{<-} "o";"3" <0pt>
\ar @{<-} "o";"4" <0pt>
\endxy}\Ea
 }_{j \times} \,.
\]
If the vertex $v$ already has an incoming edge, then there is one additional class obtained by not adding any input forest.

Hence we find that $H(\gr C_{\tilde \gamma})$ is spanned by graphs obtained from $\tilde \gamma$ as follows:
\begin{enumerate}
 \item Add some bivalent vertices with one input and one output on edges.
 \item For each vertex that is either not at least trivalent or does not have an incoming edge, sum over all ways of attaching incoming legs at that vertex.
 \item For at least trivalent vertices with an incoming edge, there is a choice of either not adding anything at that vertex, or summing over all ways of attaching incoming legs at that vertex. Let us call the vertices for which the first choice is made bald vertices and the others hairy.
\end{enumerate}

Let us look at the next page in the spectral sequence associated to our filtration on $C_{\tilde \gamma}$.
The differential creates one edge vertex by either splitting an existing edge vertex or by splitting a skeleton vertex. Again, the complex splits into a product of complexes, one for each edge of $\tilde \gamma$. For each such edge we have to consider 3 cases separately:
\begin{enumerate}
 \item Both endpoints in $\tilde \gamma$ are hairy.
 \item Both endpoints in $\tilde \gamma$ are bald.
 \item One endpoints is hairy and one is bald.
\end{enumerate}

We leave it to the reader to check that:
\begin{enumerate}
 \item In the first case the cohomology is one-dimensional, represented by a single edge without edge vertices.
 \item In the third case the cohomology vanishes.
\end{enumerate}

Since there is necessarily at least one hairy vertex in the graph, the second assertion implies that if there is a bald vertex as well, the resulting complex is acyclic.
Hence all vertices must be hairy. By the first assertion the cohomology is one-dimensional for each skeleton.
One easily checks that the representative is exactly the image of the skeleton considered as element in $\hGCor_{c+d+1}$. Hence the proposition follows.
\end{proof}

\subsubsection{\bf Remark}\label{rem:alternative Fqiso proof}
 There is also an alternative way of computing the cohomology of the deformation complex $\Der(\HoLBcd)$ (and of its completion).
 Namely, by Koszulness of $\LBcd$ this complex is quasi-isomorphic to $\Def(\HoLBcd\to \LBcd)[1]$.
 It is well known that the prop governing Lie bialgebras $\LieBiP$ may be written as
 \[
\LieBiP(n,m) \cong \bigoplus_N \caL ieP(n,N) \otimes_{\bS_N} \caL ieCP(N,m)
 \]
using the props governing Lie algebras and Lie coalgebras. Interpreting elements of the above prop as linear combinations directed acyclic graphs, the sub-properad $\LieBi$ may be obtained as that formed by the connected such graphs.
 It is hence an easy exercise to check that $\Def(\HoLBcd \to \LBcd)[1]$ is identical to the complex $\Def(\hoe_{c+d} \to e_{c+d})_{\rm conn}[1]$ from \cite{Wi1}, up to unimportant completion issues.
 The cohomology of the latter complex has been computed in loc. cit. to be
 \[
  H(\GC_{c+d}) \oplus \bigoplus_{j\geq 1\atop j\equiv 2(c+d)+1 \bmod 4} \K[c+d-j] \oplus \K \, .
 \]
Using the main result of \cite{Wi2} this agrees with the cohomology as computed by Theorem {\ref{thm:Fqiso}}. Conversely, the above proof of Theorem {\ref{thm:Fqiso}} together with this remark yields an alternative proof of the main result of \cite{Wi2}.

\subsection{The proof of Theorem {\ref{thm:Fhbarqiso}}}\label{app:defproof2}
Let us next consider Theorem {\ref{thm:Fhbarqiso}}, whose proof will be a close analog of that of Theorem {\ref{thm:Fqiso}} in the previous subsection.
There is a natural differential graded Lie algebra structure on $\hGCor_{c+d+1}$ such that the map $\Psi: \GCor_{c+d+1}\to \hGCor_{c+d+1}$ from the previous section is a map of Lie algebras. The map $\Psi$ extends $\hbar$-linearly to a map of graded Lie algebras
$$
\Psi_\hbar: \GCor_{c+d+1}[[\hbar]]\to \hGCor_{c+d+1}[[\hbar]].
$$
The Maurer-Cartan element $\Phi_\hbar\in \GCor_{c+d+1}[[\hbar]]$ from (\ref{2: Phi_hbar MC element}) is sent to a Maurer-Cartan element $\hat \Phi_\hbar:=\Psi_\hbar(\Phi_\hbar) \in \hGCor_{c+d+1}[[\hbar]]$. We endow $\hGCor_{c+d+1}[[\hbar]]$ with the differential
\[
 d_\hbar \Gamma = [\hat \Phi_\hbar, \, ].
\]
In particular it follows that we have a map of differential graded Lie algebras
\[
 \Psi_\hbar: (\GCor_{c+d+1}[[\hbar]], d_\hbar)\to (\hGCor_{c+d+1}[[\hbar]], d_\hbar).
\]

The map $F_\hbar$ from Theorem {\ref{thm:Fhbarqiso}} factors through $\hGCor_{c+d+1}[[\hbar]]$:
\[
 \GC_{c+d+1}^{or}[[\hbar]]\stackrel{\Psi_\hbar}{\longrightarrow} \hGCor_{c+d+1}[[\hbar]] \stackrel{G_\hbar}{\longrightarrow} \Der(\LoBcd).
\]
The second map $G_\hbar: \hGCor_{c+d+1}[[\hbar]]\to \Der(\HoLBcd)$ sends $\hbar^N\Gamma$, for $\Gamma\in \hGCor_{c+d+1}$ to
\[
  \sum_{j\geq 1}\sum_{k=0}^N \hbar^{N-k}
\Ba{c} \xy
(0,-5.5)*+{...},
%
(0,0)*+{_{\Ga_k}}*\cir{}="o",
(-5,-7)*{}="1",
(-3,-7)*{}="2",
(3,-7)*{}="3",
(5,-7)*{}="4",
(-3,5)*{}="5",
(3,5)*{}="6",
(5,5)*{}="7",
(-5,5)*{}="8",
\ar @{<-} "o";"1" <0pt>
\ar @{<-} "o";"2" <0pt>
\ar @{<-} "o";"3" <0pt>
\ar @{<-} "o";"4" <0pt>
\endxy
\Ea
\]
where we  again sum over all ways of attaching the incoming legs, setting to zero graphs with vertices without incoming edges. Furthermore, $\Gamma_k$ is the linear combination of graphs obtained by summing over all ways of assigning weights to the vertices of $\Gamma$ such that the total weight is $k$. 
We have the following two results, from which Theorem {\ref{thm:Fhbarqiso}} immediately follows.

\subsubsection{\bf Proposition}
{\em The map $\Psi_\hbar: (\GC_{c+d+1}^{or}[[\hbar]],d_\hbar) \to (\hGCor_{c+d+1}[[\hbar]],d_\hbar)$ is a quasi-isomorphism up to the classes $T\K[[\hbar]]\subset \hGCor_{c+d+1}[[\hbar]]$ where
 \[
  T=
  \sum_{m,p}
 (m+2p-2)
 \overbrace{
 \xy
(0,4.5)*+{...},
(0,0)*+{_p}*\cir{}="o",
(-5,-5)*{}="1",
(-3,-5)*{}="2",
(3,-5)*{}="3",
(5,-5)*{}="4",
(-3,5)*{}="5",
(3,5)*{}="6",
(5,5)*{}="7",
(-5,5)*{}="8",
%
\ar @{->} "o";"5" <0pt>
\ar @{->} "o";"6" <0pt>
\ar @{->} "o";"7" <0pt>
\ar @{->} "o";"8" <0pt>
\endxy
 }^{m\times}.
 \]
}

\begin{proof}[Proof sketch]
 Take filtrations on $\GCor_{c+d+1}[[\hbar]]$ and $\hGCor_{c+d+1}[[\hbar]]$ by the power of $\hbar$.
 The differential on the associated graded spaces is the $\hbar$-linear extension of the differentials on $\GCor_{c+d+1}$ and $\hGCor_{c+d+1}$. Hence by Proposition {\ref{prop:GChGC}} the map $\Psi_\hbar$ is a quasi-isomorphism on the level of the associated graded spaces, up to the classes above. The result follows by a standard spectral sequence argument, noting that the above element $T$ is indeed $d_\hbar$-closed.
\end{proof}

\subsubsection{\bf Proposition}\label{prop:Ghbarqiso}
{\em The map $G_\hbar: \hGCor_{c+d+1}[[\hbar]]\to \Der(\wHoLoBcd)$ is a quasi-isomorphism.
}

\begin{proof}[Proof sketch]
 Take filtrations on $\hGCor_{c+d+1}[[\hbar]]$ and
 \[
  \Der(\wHoLoBcd)\cong
  \prod_{m,n\geq 1} \left(\wLoBcd(m,n)\otimes \sgn_m^{\ot |c|}\otimes \sgn_n^{\ot |d|} \right)^{\bS_m\times \bS_n} [[\hbar]][c(1-m)+d(1-n)]
 \]
by genus and by powers of $\hbar$.
Then we claim that the induced map on the associated graded complexes $\gr G_\hbar: \gr\hGCor_{c+d+1}[[\hbar]]\to \gr\Der(\wHoLoBcd)$ is a quasi-isomorphism, thus showing the proposition by a standard spectral sequence argument.

\sip

To show the claim, we proceed analogously to the proof of Proposition {\ref{prop:Gqiso}}. Let us go through the proof again and highlight only the differences.
The skeleton of a graph is defined as before, except that one also forgets the weights of all vertices.
The complex $\gr \Der(\wHoLoBcd)$ splits into a product of subcomplexes that we again call $\tilde C_\gamma$, one for each skeleton $\gamma$. Again
\[
 \tilde C_\gamma =  C_{\tilde \gamma}^{\Aut_\gamma}
\]
for some representative $\tilde \gamma$ of the isomorphism class $\gamma$. Hence it again suffices to compute the cohomology of $C_{\tilde \gamma}$.
Graphs contributing are again obtained by adding edge vertices and input forests, except that now all vertices are also assigned an arbitrary weight.
Again we take a filtration on the number of edge vertices, which leaves us with the task of computing
the cohomology of a complex associated to one forest attached to a vertex $v$. We find that representatives of cohomology classes are either:
\begin{itemize}
 \item Vertex $v$ with any weight and no attached forest. Let us call such a $v$ again bald.
 \item Vertex $v$ with weight $0$ and input legs attached in all possible ways, let us call such a $v$ again hairy.
\end{itemize}

The differential on the second page of the spectral sequence again adds one edge vertex, which however can have a non-zero weight now, and if it has a non-zero weight it may be bald.
We may introduce another filtration by the number of non-hairy edge vertices. The differential on the associated graded creates one hairy edge vertex.
The resulting complex is a tensor product of complexes, one for each edge. The complexes associated to each edge again can have three different types: (i) both endpoints in the skeleton are hairy, (ii) both are bald or (iii) one is hairy, one is bald. Again one checks that in case (iii) the complex is acyclic and in case (i) one-dimensional, the cohomology class represented by a single edge.
Hence, since at least one vertex must be hairy, all vertices must be. Hence we recover at this stage the image of $\hGCor_{c+d+1}[[\hbar]]$ and are done.
\end{proof}

%

\subsection{Some applications}\label{4: subsec on applications} As remarked in the Introduction,
the most interesting cases for applications are ``classical" properads
$$
\LB:=\LB_{1,1}, \ \ \ \LoB:=\LoB_{1,1},  \ \ \ \ \LB_{odd}:= \LB_{0,1}.
$$
In the even case the Lie and co-Lie generators of these Lie bialgebra properads have homological degree zero, in the odd case the co-Lie generator has degree 1 and the Lie generator the degree zero.
In all three cases the associated dg Lie algebras
$$
\Der(\wHoLB)\ \  , \ \ \ \ \Der(\wHoLoB\hspace{-2mm})  \ \  , \ \ \ \  \Der(\wHoLB_{odd})
$$
are generated by graphs $\Ga$ with $\# V(\Ga) + \# E(\Ga)\geq 1$ so that all the three
are {\em positively}\, graded dg Lie algebras with respect to this parameter, and it makes sense to talk about the {\em groups}
$$
\exp\left(\Der^0(\wHoLB)\right)\ \  , \ \ \ \ \exp\left(\Der^0(\wHoLoB\hspace{-2mm})\right)  \ \  , \ \ \ \  \exp\left(\Der^0(\wHoLB_{odd})\right)
$$
which can be identified with, respectively,  degree zero subalgebras $\Der^0(\wHoLB)$, $\Der^0(\wHoLoB\hspace{-2mm})$, $\Der^0(\wHoLB_{odd})$ equipped with the standard Baker-Campbell-Hausdorff  multiplication.
The subsets of co-cycles
$$
Z^0(\wHoLB)\subset \Der^0(\wHoLB), \ \   Z^0(\wHoLoB\hspace{-2mm})\subset \Der^0(\wHoLoB\hspace{-2mm}), \ \
Z^0(\wHoLB)\subset \Der^0(\wHoLB_{odd})
$$
are precisely the automorphism groups of the completed properads,
$$
\Aut(\wHoLB)\ \  , \ \ \ \ \Aut(\wHoLoB\hspace{-2mm})  \ \  , \ \ \ \  \Aut(\wHoLB_{odd}),
$$
with the zero elements corresponding to the identity automorphisms.
We say that two automorphisms $f,g\in \Aut(...)$ from the list above are homotopy equivalent, $f\sim g$, if they differ
by a coboundary, $f-g=dh$ as elements in the corresponding dg Lie algebra $\Der(...)$, i.e.
if $f$ and $g$ define the same cohomology classes in $H^0(\Der(...))$. The subset of automorphisms homotopy equivalent to the identity automorphism (i.e.\ to zero in $\Der(...)$) is a normal subgroup
in $\Aut(...)$, and the quotient by this normal subgroup is called the {\em group of homotopy non-trivial}\, automorphisms and is  denoted, respectively, by
\Beq\label{4: three auto groups}
\mathbf{Aut}(\wHoLB)\ \  , \ \ \ \ {\mathbf{Aut}}(\wHoLoB\hspace{-2mm})  \ \  , \ \ \ \  \mathbf{Aut}(\wHoLB_{odd}).
\Eeq
Then we have the following three corollaries to the Main Theorems.

\subsubsection{\bf Proposition} (i) {\em The group $\mathbf{Aut}(\wHoLB)$ is equal to the Grothendieck-Teichm\"uller group $GRT=GRT_1\ltimes \K^*$ with the subgroup $\K^*$
 acting on $\wHoLB$ by the following rescaling transformations of the generators},
 \Beq\label{4: rescalings of the generators}
 \Ba{c}\resizebox{12mm}{!}{\begin{xy}
 <0mm,0mm>*{\circ};<0mm,0mm>*{}**@{},
 <-0.6mm,0.44mm>*{};<-8mm,5mm>*{}**@{-},
 <-0.4mm,0.7mm>*{};<-4.5mm,5mm>*{}**@{-},
 <0mm,0mm>*{};<-1mm,5mm>*{\ldots}**@{},
 <0.4mm,0.7mm>*{};<4.5mm,5mm>*{}**@{-},
 <0.6mm,0.44mm>*{};<8mm,5mm>*{}**@{-},
   <0mm,0mm>*{};<-8.5mm,5.5mm>*{^1}**@{},
   <0mm,0mm>*{};<-5mm,5.5mm>*{^2}**@{},
   <0mm,0mm>*{};<4.5mm,5.5mm>*{^{m\hspace{-0.5mm}-\hspace{-0.5mm}1}}**@{},
   <0mm,0mm>*{};<9.0mm,5.5mm>*{^m}**@{},
 <-0.6mm,-0.44mm>*{};<-8mm,-5mm>*{}**@{-},
 <-0.4mm,-0.7mm>*{};<-4.5mm,-5mm>*{}**@{-},
 <0mm,0mm>*{};<-1mm,-5mm>*{\ldots}**@{},
 <0.4mm,-0.7mm>*{};<4.5mm,-5mm>*{}**@{-},
 <0.6mm,-0.44mm>*{};<8mm,-5mm>*{}**@{-},
   <0mm,0mm>*{};<-8.5mm,-6.9mm>*{^1}**@{},
   <0mm,0mm>*{};<-5mm,-6.9mm>*{^2}**@{},
   <0mm,0mm>*{};<4.5mm,-6.9mm>*{^{n\hspace{-0.5mm}-\hspace{-0.5mm}1}}**@{},
   <0mm,0mm>*{};<9.0mm,-6.9mm>*{^n}**@{},
 \end{xy}}\Ea
  \lon \la^{m+n-2}
 \Ba{c}\resizebox{12mm}{!}{\begin{xy}
 <0mm,0mm>*{\circ};<0mm,0mm>*{}**@{},
 <-0.6mm,0.44mm>*{};<-8mm,5mm>*{}**@{-},
 <-0.4mm,0.7mm>*{};<-4.5mm,5mm>*{}**@{-},
 <0mm,0mm>*{};<-1mm,5mm>*{\ldots}**@{},
 <0.4mm,0.7mm>*{};<4.5mm,5mm>*{}**@{-},
 <0.6mm,0.44mm>*{};<8mm,5mm>*{}**@{-},
   <0mm,0mm>*{};<-8.5mm,5.5mm>*{^1}**@{},
   <0mm,0mm>*{};<-5mm,5.5mm>*{^2}**@{},
   <0mm,0mm>*{};<4.5mm,5.5mm>*{^{m\hspace{-0.5mm}-\hspace{-0.5mm}1}}**@{},
   <0mm,0mm>*{};<9.0mm,5.5mm>*{^m}**@{},
 <-0.6mm,-0.44mm>*{};<-8mm,-5mm>*{}**@{-},
 <-0.4mm,-0.7mm>*{};<-4.5mm,-5mm>*{}**@{-},
 <0mm,0mm>*{};<-1mm,-5mm>*{\ldots}**@{},
 <0.4mm,-0.7mm>*{};<4.5mm,-5mm>*{}**@{-},
 <0.6mm,-0.44mm>*{};<8mm,-5mm>*{}**@{-},
   <0mm,0mm>*{};<-8.5mm,-6.9mm>*{^1}**@{},
   <0mm,0mm>*{};<-5mm,-6.9mm>*{^2}**@{},
   <0mm,0mm>*{};<4.5mm,-6.9mm>*{^{n\hspace{-0.5mm}-\hspace{-0.5mm}1}}**@{},
   <0mm,0mm>*{};<9.0mm,-6.9mm>*{^n}**@{},
 \end{xy}}\Ea
 \ \ \ \ \ \ \ \ \ \forall\ \la\in \K^*.
 \Eeq

(ii) {\em The group  ${\mathbf{Aut}}(\wHoLoB\hspace{-2mm})$ is equal to $GRT=GRT_1\ltimes \K^*$ with the subgroup $\K^*$ acting on $\wHoLB$ by  rescaling transformations},
$$
\Ba{c}\resizebox{12mm}{!}{\xy
(-9,-6)*{};
(0,0)*+{a}*\cir{}
**\dir{-};
(-5,-6)*{};
(0,0)*+{a}*\cir{}
**\dir{-};
(9,-6)*{};
(0,0)*+{a}*\cir{}
**\dir{-};
(5,-6)*{};
(0,0)*+{a}*\cir{}
**\dir{-};
(0,-6)*{\ldots};
(-10,-8)*{_1};
(-6,-8)*{_2};
(10,-8)*{_n};
(-9,6)*{};
(0,0)*+{a}*\cir{}
**\dir{-};
(-5,6)*{};
(0,0)*+{a}*\cir{}
**\dir{-};
(9,6)*{};
(0,0)*+{a}*\cir{}
**\dir{-};
(5,6)*{};
(0,0)*+{a}*\cir{}
**\dir{-};
(0,6)*{\ldots};
(-10,8)*{_1};
(-6,8)*{_2};
(10,8)*{_m};
\endxy}\Ea
\lon
\la^{m+n+a-2}
\Ba{c}\resizebox{12mm}{!}{\xy
(-9,-6)*{};
(0,0)*+{a}*\cir{}
**\dir{-};
(-5,-6)*{};
(0,0)*+{a}*\cir{}
**\dir{-};
(9,-6)*{};
(0,0)*+{a}*\cir{}
**\dir{-};
(5,-6)*{};
(0,0)*+{a}*\cir{}
**\dir{-};
(0,-6)*{\ldots};
(-10,-8)*{_1};
(-6,-8)*{_2};
(10,-8)*{_n};
(-9,6)*{};
(0,0)*+{a}*\cir{}
**\dir{-};
(-5,6)*{};
(0,0)*+{a}*\cir{}
**\dir{-};
(9,6)*{};
(0,0)*+{a}*\cir{}
**\dir{-};
(5,6)*{};
(0,0)*+{a}*\cir{}
**\dir{-};
(0,6)*{\ldots};
(-10,8)*{_1};
(-6,8)*{_2};
(10,8)*{_m};
\endxy}\Ea
$$

(iii) {\em The group  ${\mathbf{Aut}}(\wHoLB_{odd})$ is equal to $\K^*$ which acts $\wHoLB_{odd}$ by rescaling transformations of the generators as in (\ref{4: rescalings of the generators})}.

\begin{proof} All three groups in (\ref{4: three auto groups}) can be identified as sets with, respectively, the zero-th  cohomology groups
$$
H^0(\Der(\wHoLB))\ \  , \ \ \ \ H^0(\Der(\wHoLoB\hspace{-2mm}))  \ \  , \ \ \ \  H^0(\Der(\wHoLB_{odd}))
$$
which in turn, by the Main Theorems, can be identified with the zero-th cohomology groups
$$
H^0(\GCor_3, \delta) \oplus \K \ , \ \ \ H^0(\GCor_3[[\hbar]], \delta_\hbar)\oplus \K \ , \ \ \ H^0(\GCor_2, \delta)\oplus \K,
$$
which in turn, according to \cite{Wi1}, \cite{CMW} and, respectively, \cite{Wi2}, are equal as Lie algebras to
$$
\fg\fr\ft_1 \oplus \K, \ \ \ \ \ \  \ \fg\fr\ft_1 \oplus \K,\ \ \ \mbox{{and\ respectively}} \ \  \ \ \K.
$$
This proves all the claims.
\end{proof}

This Proposition implies a highly non-trivial action of the Grothendieck-teichm\"uller group $GRT_1$
on the completed properads $\wHoLB$ and $\wHoLoB$\hspace{-1mm} and hence on their representations.
However one must be careful when talking about representations of the {\em completed}\, free properads. We introduce a topology on properads $\wHoLB/\wHoLB_{odd}$ (respectively on $\wHoLoB$\hspace{-2mm}) as the one induced by the genus filtration (respectively, by the filtration defined by the parameter ``genus + total $a$-weight", cf. \cite{CMW}). If $W$ is a topological vector space, then by a representation of $\wHoLB/\wHoLB_{odd}$ (respectively, $\wHoLoB$\hspace{-2mm}) in $W$ we mean a {\em continuous}\, map of topological properads
$$
\rho: \wHoLB/\wHoLB_{odd} \lon \cE nd_W, \ \ \ \ \ \rho^\diamond:\wHoLoB \lon \cE nd_W.
$$
Let us describe a sufficiently large class of such representations. Consider an arbitrary dg space  $V$, and let $\hbar$ be formal parameter of degree zero. The vector space
$W:=V[[\hbar]]$ can be equipped with the standard $\hbar$-adic topology. Then maps $\rho$ and $\rho^\diamond$ satisfying on the generators the condition
$$
\rho\left(
\Ba{c}\resizebox{12mm}{!}{\begin{xy}
 <0mm,0mm>*{\circ};<0mm,0mm>*{}**@{},
 <-0.6mm,0.44mm>*{};<-8mm,5mm>*{}**@{-},
 <-0.4mm,0.7mm>*{};<-4.5mm,5mm>*{}**@{-},
 <0mm,0mm>*{};<-1mm,5mm>*{\ldots}**@{},
 <0.4mm,0.7mm>*{};<4.5mm,5mm>*{}**@{-},
 <0.6mm,0.44mm>*{};<8mm,5mm>*{}**@{-},
   <0mm,0mm>*{};<-8.5mm,5.5mm>*{^1}**@{},
   <0mm,0mm>*{};<-5mm,5.5mm>*{^2}**@{},
   <0mm,0mm>*{};<4.5mm,5.5mm>*{^{m\hspace{-0.5mm}-\hspace{-0.5mm}1}}**@{},
   <0mm,0mm>*{};<9.0mm,5.5mm>*{^m}**@{},
 <-0.6mm,-0.44mm>*{};<-8mm,-5mm>*{}**@{-},
 <-0.4mm,-0.7mm>*{};<-4.5mm,-5mm>*{}**@{-},
 <0mm,0mm>*{};<-1mm,-5mm>*{\ldots}**@{},
 <0.4mm,-0.7mm>*{};<4.5mm,-5mm>*{}**@{-},
 <0.6mm,-0.44mm>*{};<8mm,-5mm>*{}**@{-},
   <0mm,0mm>*{};<-8.5mm,-6.9mm>*{^1}**@{},
   <0mm,0mm>*{};<-5mm,-6.9mm>*{^2}**@{},
   <0mm,0mm>*{};<4.5mm,-6.9mm>*{^{n\hspace{-0.5mm}-\hspace{-0.5mm}1}}**@{},
   <0mm,0mm>*{};<9.0mm,-6.9mm>*{^n}**@{},
 \end{xy}}\Ea
 \right)\in \hbar^{m+n-2}\Hom(V^{\ot n}, V^{\ot m}))[[\hbar]], \ \ \ \ \ \ \ \
 \rho^\diamond\left(
 \Ba{c}\resizebox{12mm}{!}{\xy
(-9,-6)*{};
(0,0)*+{a}*\cir{}
**\dir{-};
(-5,-6)*{};
(0,0)*+{a}*\cir{}
**\dir{-};
(9,-6)*{};
(0,0)*+{a}*\cir{}
**\dir{-};
(5,-6)*{};
(0,0)*+{a}*\cir{}
**\dir{-};
(0,-6)*{\ldots};
(-10,-8)*{_1};
(-6,-8)*{_2};
(10,-8)*{_n};
(-9,6)*{};
(0,0)*+{a}*\cir{}
**\dir{-};
(-5,6)*{};
(0,0)*+{a}*\cir{}
**\dir{-};
(9,6)*{};
(0,0)*+{a}*\cir{}
**\dir{-};
(5,6)*{};
(0,0)*+{a}*\cir{}
**\dir{-};
(0,6)*{\ldots};
(-10,8)*{_1};
(-6,8)*{_2};
(10,8)*{_m};
\endxy}\Ea
\right)\in \hbar^{m+n+a-2}\Hom(V^{\ot n}, V^{\ot m}))[[\hbar]]
$$
define a\, {\em continuous}\, representation of the properads $\wHoLB/\wHoLB_{odd}$ and, respectively, $\wHoLoB$\hspace{-2mm} in the topological space $V[[\hbar]]$.

\sip

\subsubsection{\bf On the unique non-trivial deformation of $\wHoLB_{odd}$}
It was proven in \cite{Wi2} that the cohomology group $H^1(\GCor_2)$ is one-dimensional and is spanned by the following graph
$$
\Upsilon_4:=\Ba{c}\xy
(0,0)*{\bullet}="1",
(-7,16)*{\bullet}="2",
(7,16)*{\bullet}="3",
(0,10)*{\bullet}="4",
\ar @{<-} "2";"4" <0pt>
\ar @{<-} "3";"4" <0pt>
\ar @{<-} "4";"1" <0pt>
\ar @{<-} "2";"1" <0pt>
\ar @{<-} "3";"1" <0pt>
\endxy\Ea
+
2
\Ba{c}\xy
(0,0)*{\bullet}="1",
(-6,6)*{\bullet}="2",
(6,10)*{\bullet}="3",
(0,16)*{\bullet}="4",
\ar @{<-} "4";"3" <0pt>
\ar @{<-} "4";"2" <0pt>
\ar @{<-} "3";"2" <0pt>
\ar @{<-} "2";"1" <0pt>
\ar @{<-} "3";"1" <0pt>
\endxy\Ea
+
 \Ba{c}\xy
(0,16)*{\bullet}="1",
(-7,0)*{\bullet}="2",
(7,0)*{\bullet}="3",
(0,6)*{\bullet}="4",
\ar @{->} "2";"4" <0pt>
\ar @{->} "3";"4" <0pt>
\ar @{->} "4";"1" <0pt>
\ar @{->} "2";"1" <0pt>
\ar @{->} "3";"1" <0pt>
\endxy\Ea.
$$

Moreover $H^2(\GCor_2)=\K$ and is spanned by a graph with four vertices. This means that one can construct by induction a new Maurer-Cartan element in the Lie algebra $\GCor_2$ (the integer subscript in the summand $\Upsilon_n$ stands for the number of vertices of graphs)
$$
\Upsilon_{KS}= \xy
 (0,0)*{\bullet}="a",
(6,0)*{\bu}="b",
\ar @{->} "a";"b" <0pt>
\endxy  + \Upsilon_4
+ \Upsilon_6 + \Upsilon_8 + \ldots
$$
as all obstructions have more than $7$ vertices and hence do not hit the unique cohomology class
in $H^2(\GCor_2)$. Up to gauge equivalence, this new Maurer-Cartan element $\Upsilon$ is the {\em only}\, non-trivial deformation of the standard Maurer-Cartan element $\xy
 (0,0)*{\bullet}="a",
(6,0)*{\bu}="b",
\ar @{->} "a";"b" <0pt>
\endxy$.
We call this element  {\em Kontsevich-Shoikhet}\, one as it was first found by Boris Shoikhet in \cite{Sh} with a reference to an important contribution by Maxim Kontsevich via an informal communication.

\sip

By Main theorem {\ref{thm:Fqiso}}, the Maurer-Cartan element $\Upsilon_{KS}$ equips the completed
non-differential properad $\wHoLB_{odd}$  with a new differential denoted by
$\delta_{KS}$. If continuous representations of $\wHoLB_{odd}$ equipped with the standard differential $\delta$ (originating from $\xy
 (0,0)*{\bullet}="a",
(6,0)*{\bu}="b",
\ar @{->} "a";"b" <0pt>
\endxy$) in a topological vector space $V[[\hbar]]$ can be identified with ordinary formal Poisson structures $\pi\in \cT_{poly}(V)[[\hbar]]$, the continuous representations of $\wHoLB_{odd}$ equipped with the new differential $\delta_{KS}$  give us a notion of {\em quantizable Poisson structure}\, $\pi^{quant}\in \cT_{poly}(V)[[\hbar]]$ (this notion can be globalized from a vector space $V$ to an arbitrary manifold $M$). It was proven in \cite{MW3} that for {\em finite}-dimensional vector spaces
$V$ (or manifolds $M$), there is a one-to-one correspondence between ordinary Poisson structures and quantizable ones, but this correspondence
$$
\left\{ \mbox{Ordinary Poisson structures}\ \pi\ \mbox{on}\ M \right\}
\stackrel{1:1}{\leftrightarrow}
\left\{ \mbox{Quantizable Poisson structures}\ \pi^{quant} \ \mbox{on}\ M\right\}
$$
is highly non-trivial and depends on a choice of Drinfeld associator. Moreover, quantizable
Poisson structures can be deformation quantized in a trivial (and essentially unique) perturbative way \cite{MW3} so that all the subtleties of the deformation quantization are hidden  in the above correspondence.

\def\cprime{$'$}

\end{document}